\documentclass[11pt]{amsart}
\usepackage{amssymb,amscd,amsthm,amsxtra}
\usepackage{latexsym}

\newtheorem{thm}{Theorem}[section]

\newtheorem{lem}[thm]{Lemma}
\newtheorem{prop}[thm]{Proposition}
\theoremstyle{definition}

\theoremstyle{remark}
\newtheorem{rem}[thm]{Remark}
\numberwithin{equation}{section}

\newcommand{\comment}[1]{}

\def\a{\alpha}
\def\q{\frac{1}{2}}

\def\R{\mathbb{R}}
\def\C{\mathbb{C}}
\def\Z{\mathbb{Z}}

\def\beq{\begin{equation}}
\def\eeq{\end{equation}}

\def\beq{\begin{equation}}
\def\eeq{\end{equation}}
\def\ra{\rangle}
\def\la{\langle}
\def\ri{\rightarrow}
\def\les{\lesssim}

\begin{document}

\title[Low regularity solutions for a 2D quadratic NLS]
{Low regularity solutions for a 2D quadratic non-linear
Schr\"odinger equation}
\author{Ioan Bejenaru}
\author{Daniela De Silva}



\begin{abstract}
We establish that the initial value problem for the quadratic
non-linear Schr\"odinger equation
$$ iu_t - \Delta u = u^2$$
where $u: \R^2 \times \R \to \C$, is locally well-posed in
$H^s(\R^2)$ when $s > -1$. The critical exponent for this problem
is $s_c=-1$ and previous work in \cite{c1} established local
well-posedness for $s > -3/4$.
\end{abstract}
\maketitle \footnotetext[1]{Both authors were partial supported by
the Mathematical Sciences Research Institute (MSRI) in Berkeley. }
\section{Introduction}

The aim of this paper is to continue the development of the theory
for quadratic non-linear Schr\"odinger equations of the form:
\begin{equation} \label{E}
\begin{cases}
\begin{aligned}
iu_{t}-\Delta u &= P(u,\bar{u}), \   t\in \mathbb{R}, x\in \mathbb{R}^{n} \\
u(x,0) &=u_{0}(x) \in H^s(\mathbb{R}^n)
\end{aligned}
\end{cases}
\end{equation}
where $u: \mathbb{R}^{n} \times \mathbb{R} \rightarrow
\mathbb{C},$ and $P: \mathbb{C}^{2} \rightarrow \mathbb{C}$ is a
quadratic polynomial. Here $H^s(\mathbb{R}^n)$ denotes the usual
inhomogeneous Sobolev space.

An important concept for this type of problem is the scaling
(critical) exponent $s_{c}$. This is the exponent of the Sobolev
space which scales in the same way as the equation. A
straightforward computation shows that the critical exponent for
problem (\ref{E}) is $s_{c}=\frac{n}{2}-2$. Heuristically one
would expect to obtain a local well-posedness theory for initial
data $u_0 \in H^s$ for all $s \geq s_{c}$. For a precise
definition of local well-posedness, we refer to the statement of
our main result, Theorem \ref{lwp}.

As for many other evolution equations of this type, it turns out
that in lower dimensions obtaining local well-posedness for $s
\geq s_c$ is a delicate problem. Indeed, for $n \geq 4$ it has
been shown in \cite{cwi}, via Strichartz estimates, that local
well-posedness holds for (\ref{E}) for all $s \geq s_c$, while for
$n \leq 3$ local well-posedness holds for all $s \geq 0$. On the
other hand, for $n \leq 3$, $s_c =\frac{n}{2}-2 < 0$, therefore it
is expected that one needs more refined techniques in order to get
closer to the scaling exponent.

One new feature in low dimension is that one has to distinguish
among the three types of nonlinearities $u^2$, $|u|^2=u \bar{u}$
and $\bar{u}^2$. These nonlinearities behave differently and the
local well-posedness results obtained for each of them are not the
same.

In the case $n=3$ ($s_c=-\q$) \cite{tao1} established local
well-posedness for $s > -\q$ for the non-linearities $u^2$ and
$\bar{u}^2$, while local well-posedness for the non-linearity
$|u|^2$ is shown to hold for $s
> -\frac{1}{4}$. The open problem in $3D$ is what
happens at scaling, i.e. $s=-\q$, for $u^{2}$ and $\bar{u}^2$ and
at $s=-\frac{1}{4}$ for $|u|^2$. We have strong reasons to believe
that for the nonlinearity $|u|^2$ one cannot get a positive result
for $s < -\frac{1}{4}$.

In dimensions one and two, the available results are even further
away from scaling. Precisely, in one dimension
($s_c=-\frac{3}{2}$) the following results are known:

\begin{itemize}

\item for $u^{2}$: well-posedness for $s \geq -1$ and ill-posedness for $s<-1$, see
\cite{bt};

\item for $|u|^2$: well-posedness for $s > -\frac{1}{4}$, see
\cite{kpv};

\item for $\bar{u}^2$: well-posedness for $s > -\frac{3}{4}$, see
\cite{kpv}.

\end{itemize}

We remark that the technique in \cite{bt} is not directly
applicable to improve the results known for the nonlinearities
$|u|^2$ and $\bar{u}^2$. Indeed such technique is based on
``undoing" the counterexample in \cite{kpv} imposing the condition
$s>-3/4$ for the nonlinearity $u^2$. For $|u|^2$, the authors of
\cite{bt} could not undo the counterexample yielding $s
> -\frac{1}{4}$, while the nonlinearity $\bar{u}^2$ is
known as being the easiest of all and a simpler argument should
suffice. Hence, work is needed for these two cases.

For $n=2$ the state of art was established in \cite{c1}. For the
nonlinearities $u^{2}$ and $\bar{u}^2$ the authors proved
well-posedness for $s > -\frac{3}{4}$, while for the nonlinearity
$|u|^2$ well-posedness is shown to hold for $s > -\frac{1}{4}$.

Before the result in \cite{bt}, the breakpoints for $u^2$ in $1D$
and $2D$, namely $s> -\frac{3}{4}$, were imposed by a very similar
counterexample, see \cite{kpv} and \cite{c1} respectively.
Therefore we decided to investigate a possible improvement for the
nonlinearity $u^2$ in $2D$ following the ideas in \cite{bt}. We
have used similar structures for the functional spaces introduced
to overcome the deadlock imposed by the known counterexample in
\cite{c1}. The techniques are more involved since we deal with the
two dimensional problem which brings additional difficulties.

Recall that we are interested in the problem:
\begin{equation} \label{EE}
\begin{cases}
\begin{aligned}
iu_{t}-\Delta u &= u^2, \   t\in \mathbb{R}, x\in \mathbb{R}^{2} \\
u(x,0) &=u_{0}(x) \in H^s(\mathbb{R}^2).
\end{aligned}
\end{cases}
\end{equation}
We prove the following result.

\begin{thm} \label{lwp}

Let $s>-1$ and $r > 0$ be any radius, and let $B_r$ be the ball
$$ B_r := B_{H^{s}_x(\R)}(0,r) := \{ u_0 \in H^{s}_x(\R): \| u_0 \|_{H^{s}_x(\R)} < r \}.$$
Then there exists a time $T > 0$ and a map $f \mapsto u[f]$ which
is continuous from $B_r$ to $C^0_t H^{s}_x([0,T] \times \R)$, such
that the restriction of this map to $B_r \cap H^{s'}_x(\R)$ (with
the $H^{s'}_x(\R)$ topology) maps continuously to $C^0_t
H^{s'}([0,T] \times \R)$ for any $s' \geq s$.  Furthermore, if $f$
lies in a smooth space, say $B_r \cap H^3_x(\R)$, then $u[f]$ lies
in $C^0_t H^3_x \cap C^1_t H^1_x([0,T] \times \R)$ and solves the
equation (\ref{EE}) in the classical sense.

\end{thm}

Our result does not say or predict anything about the behavior of
(\ref{EE}) at scaling, i.e. for $u_0 \in H^{-1}$. Potentially we
could prove some logarithmic divergence in the bilinear estimate
in the spirit of \cite{nak}. However this would not say almost
anything about the well/ill-posedness of the problem in this case.

Concerning the nonlinearities $|u|^2$ and $\bar{u}^2$ we believe
that our paper together with \cite{bt} provide most of the tools
needed for expanding the known results both in the $1D$ and $2D$
cases. This would establish a satisfactory well-posedness theory
in $1D$ and $2D$ for (\ref{E}).

The formalism needed for our result was developed in \cite{bt}. In
the next section, Section 2, we provide a concise but rigorous
definition of the functional space in which we need to work. For
more details about the motivation for the particular structure of
such a function space, we refer the reader to \cite{bt}. Sections
3 through 6 are devoted to the main estimates of the paper, namely
the bilinear estimates in the functional space introduced in
Section 2.

\vspace{1mm}

\subsection*{Acknowledgments.} The authors would like to thank James Colliander for
encouragement with this project.

\section{Description of the Function Space $W^s$}

In this section we introduce some notation and we define the
function space $W^s$ whose properties will allow us to develop a
local well-posedness theory for the initial value problem
\eqref{EE}, for all $s>-1.$ We fix $T=1$ and build our space
adapted to this scale. A standard rescaling argument gives us the
result for all $T$'s.

Throughout this paper we use the notation $A \les B$ to mean that
$A \leq C B$ for some $C$ which may change from line to line, but
which is independent of any of the possible variables in our
problem. If $A \lesssim B$ and $B \lesssim A$ we say that $A
\approx B$. In addition $\langle a \rangle:=1+|a|$ and $a+:=a +
\epsilon$ with $\epsilon $ positive and small.

For any $s,b\in \R$, we define $\hat X^{s,b}$ to be the closure of
the smooth functions $f : \R^{2} \times \R \rightarrow \C$ under
the following norm:
$$ \|f\|_{\hat X^{s,b}} := \| \langle \xi \rangle^s \langle \tau - |\xi|^2 \rangle^b f\|_{L^2_\xi L^2_\tau}.$$
These are the Fourier transforms of the usual $X^{s,b}$ spaces
defined in \cite{borg:xsb}.

The function space $W^s$ will be a modification of $\hat
X^{s,1/2+}$. As remarked in the introduction, a close analysis of
\cite{c1} reveals that in order to obtain well-posedness below
$s=-3/4$ one has to deal with the same kind of ``bad interactions"
which appear in the 1-dimensional case. For this reason we follow
closely \cite{bt} toward the construction of the space $W^s$.

Since $W^s$ will be constructed directly on the Fourier side (like
$\hat{X}^{s,b}$), we require that $W^s$ satisfies the following
properties.

\

$\bullet$ \it Monotonicity: \rm If $|f| \leq |g|$ pointwise, then
$\|f\|_{W^s} \leq \|g\|_{W^s}$. In particular, $\|f\|_{W^s} =
\||f|\|_{W^s}$.

\

$\bullet$ \it $H^{s}$ Energy estimate: \rm
\begin{equation} \label{es}
\| \langle \xi \rangle^{s} f \|_{L^2_\xi L^1_\tau} \les
\|f\|_{W^s}.
\end{equation}

\

$\bullet$ \it Homogeneous $H^{s}$ solution estimate: \rm
\begin{equation} \label{he}
\| f \|_{W^s} \les \| f \|_{\hat X^{s,100}}.
\end{equation}

\

$\bullet$ \it Bilinear estimate: \rm
\begin{equation} \label{be}
\| \langle \tau - |\xi|^2 \rangle^{-1} f * g \|_{W^s} \les \| f
\|_{W^s} \|g\|_{W^s},
\end{equation}
where $f*g$ denotes space-time convolution
$$ f*g(\xi,\tau):= \int_{\R} \int_{\R^2} f(\xi_1,\tau_1) g(\xi_2,\tau_2)\ d\xi_1 d\tau_1$$
using the convention
\begin{equation*}\label{tau-convention}
(\xi_1,\tau_1) + (\xi_2,\tau_2) = (\xi,\tau).
\end{equation*}

\

Once we determine a function space $W^s$ that satisfies all of the
properties above, then the machinery developed in \cite{bt} will
give the result in Theorem \ref{lwp}.

This section provides the construction of $W^s$ and the (almost)
trivial check that $W^s$ satisfies the first three properties
above. The next sections will deal with the bilinear estimate
\eqref{be}, which is indeed the core of our paper.

In order to define $W^s$, we shall partition the frequency space
$(\xi,\tau)$ in the sets $A_j \cap B_d$ for $j,d \geq 0$, where
the $A_j$'s are the annuli
$$ A_j := \{ (\xi,\tau) \in \R \times \R^2: 2^j \leq \langle \xi \rangle < 2^{j+1}\},$$
while the $B_d$'s are the parabolic neighborhoods
$$ B_d := \{ (\xi,\tau) \in \R \times \R^2: 2^d \leq \langle \tau - |\xi|^2 \rangle < 2^{d+1} \}.$$
In the next sections, we will also use the annuli $C_m$'s, $m \geq
0,$ defined by
$$ C_m := \{ (\xi,\tau) \in \R \times \R^2: 2^{m} \leq \langle \tau \rangle < 2^{m+1}\}.$$
Moreover, we denote by
$$ A_{\leq j} := \bigcup_{j' \leq j} A_{j'}; \quad B_{\leq d} := \bigcup_{d' \leq d}
B_{d'}; $$ and similarly one can define $A_{\geq j}$, $A_{>j}$,
$B_{\geq d}$, $B_{>d}$, etc.  Finally for any smooth $f$, we
denote by
$$f_j:=\chi_{A_j}f,$$
$$f_{j,d}:=\chi_{A_j\cap B_d}f,$$ where $\chi_A$ is the
characteristic function of a set $A$. We then have\footnote{All
sums and unions involving $j$ and $d$ shall be over the
non-negative integers unless otherwise mentioned.}
\begin{equation}\label{wt}
 \| f\|_{\hat X^{s,b}} \approx
(\sum_{j} \sum_{d} 2^{2sj} 2^{2bd} \| f_{j,d} \|_{L^2_\xi
L^2_\tau}^2)^{1/2}.
\end{equation}

Now define the space $\hat X^{s,b,1}$, that is the Besov
endpoint refinement of $\hat X^{s,b}$, by the following norm
\begin{equation}\label{fx}
 \| f \|_{\hat X^{s,b,1}} :=  (\sum_{j} 2^{2sj} (\sum_{d} 2^{bd} \| f_{j,d} \|_{L^2_\xi L^2_\tau})^2)^{1/2}.
\end{equation}
If $b=\q$, a straightforward computation shows that
$\hat{X}^{s,\q} \subset \hat{X}^{s,\q,1} \subset \hat{X}^{s,\q+}.$
This Besov-type space allows us to handle the ``parallel
interaction'' case, that is the case when the nonlinearity
interacts two components of the solution $u$ with the same high
frequency. However, this space alone is not sufficient even to
handle the endpoint $s=-3/4$, because of a divergence near the "$\tau$"
axis. Thus we need to further modify this space toward the
definition of $W^s$. We introduce the following function space
$Y^s$, defined via the norm
\begin{align}\label{fy}
\|f\|_{Y^s} & := \| \langle \xi \rangle^{s} f \|_{L^2_\xi
L^1_\tau}
+ \| \langle (\xi,\tau) \rangle^{s+1}f \|_{L^2_\xi L^2_\tau} \\
& :=\| f \|_{\langle \xi \rangle^{s} L^2_\xi L^1_\tau} + \| f
\|_{\langle (\xi,\tau) \rangle^{s+1} L^2_\xi L^2_\tau}\nonumber
\end{align}
where we denote by $|(\xi,\tau)|= (|\tau|+|\xi|^2)^{\frac{1}{2}}.$
Then we define the sum $$Z^s := \hat X^{s,1/2,1}+Y^s$$ which is
endowed with the usual norm
$$ \| f \|_{Z^s} := \inf \{ \| f_1 \|_{\hat X^{s,1/2,1}} + \|f_2\|_{Y^s}: f_1 \in \hat X^{s,1/2,1}; f_2 \in Y^s
; f=f_1+f_2 \}.$$ It is easy to verify that $Z^s$ is a Banach
space, with the Schwartz functions being dense. Using this space
we will be able to handle the divergences occurring along the time
axis. By definition,
$$ \|f\|_{Z^s} \leq \|f\|_{\hat X^{s,1/2,1}}, \|f\|_{Y^s}.$$
Vice-versa, in order to prove a linear estimate of the form
$\|Tf\|_{Z^s} \les \|f\|_{Z^s}$, it suffices to prove both the
estimates $\|Tf\|_{Z^s} \les\|f\|_{\hat X^{s,1/2,1}}$, and
$\|Tf\|_{Z^s} \les \|f\|_{Y^s}$.

It is a simple exercise to establish that if $f$ is smooth then
\beq \label{mm} || f\|_{\langle \xi \rangle^{s} L^2_\xi L^1_\tau}
\les \|f\|_{Z^s} \ \ \mbox{and} \ \ || f \|_{Z^s} \les \| f
\|_{\hat X^{s,90}}. \eeq

We remark that because of the $L^2$ structure of the spaces
involved in our analysis we have the following localization
property
\begin{equation}\label{localization}
\|f\|_{\mathcal{X}} \approx (\sum_j\|f_j
\|_{\mathcal{X}}^2)^{1/2},
\end{equation} for $\mathcal{X}$ either of $\hat{X}^{s,1/2,1},Y^s,$ or $Z^s$.

The two spaces $\hat X^{s,1/2,1}$ and $Y^s$ paste together nicely
along the fuzzy boundary $\langle \tau - |\xi|^2 \rangle \approx
\langle \xi \rangle^2$. More precisely, the following Lemma holds
(for details of the proof see \cite{bt}).

\begin{lem}[Pasting lemma]\label{paste}  Let $f$ be a smooth function and let $-1\leq s<0$.
If $f$ is supported on $\bigcup_{j} A_j \cap B_{\geq 2j - 100}$,
then
\begin{equation}\label{paste-a}
\|f\|_{Y^s} \les \|f\|_{Z^s}.
\end{equation}
Conversely, if $f$ is supported on $\bigcup_{j} A_j \cap B_{\leq
2j + 100}$, then
\begin{equation}\label{paste-b}
 \|f\|_{\hat X^{s,1/2,1}} \les \|f\|_{Z^s}.
 \end{equation}
\end{lem}

Now, set $$K=\cup_{j} (A_{j} \cap B_{\leq 2j-4}).$$ Notice that on
$K \cap A_{j}$ we have $2^{2j-4} \leq \tau \leq 2^{2j+4}$. For a
smooth function $f$, we define $$f^K:=\chi_{K} f, \;
f^{K^c}:=\chi_{K^c} f.$$ Here $K^c$ denotes the complement of $K.$
Then, the pasting Lemma implies that \beq \label{pp} ||f||_{Z^s}
\approx ||f^{K}||_{\hat{X}^{s,\q,1}}+||f^{K^c}||_{Y^s} \eeq This
is the way we essentially think of measuring functions in $Z^s$.

The space $Z^s$ is a good candidate for $W^s$, as it is able to
cope with two of the dangerous quadratic interactions in the
equation (namely the parallel interactions, and the interactions
which output near the time axis). There is a third type of
interaction which could cause trouble, when a solution component
near the parabola $\{ \tau = |\xi|^2 \}$ interacts with a
component near the reflected parabola $\{ \tau = -|\xi|^2 \}$ to
create a large contribution near the frequency origin. However, we do expect 
the solution to stay concentrated in the
upper half-plane $\tau
> 0$. To exploit this fact we shall introduce a weight
\begin{equation}\label{weight-def}
w(\xi,\tau) := \max(1,-\tau)^{10}
\end{equation}
to localize to the upper half-plane, and define $W^s$ to be the
space
\begin{equation}\label{W-def}
 \|f\|_{W^s} := \| w f \|_{Z^s}.
\end{equation}

The first three properties that $W^s$ is required to satisfy (see
the beginning of this section) are straightforward. The
monotonicity of $\hat X^{s,1/2,1}$,$Y^s$ and hence of $W^s$ is
immediate. The $H^s$ energy estimate (\ref{es}) follows from
(\ref{mm}) (since $w \geq 1$), while the homogeneous $H^s$
solution estimate (\ref{he}) follows directly from the following
inequalities,
$$ \|f\|_{W^s} \leq \|w f\|_{\hat X^{s,1/2,1}} \les \| w f\|_{\hat X^{s,90}} \les \|f\|_{\hat X^{s,100}},$$
where we have used the crude estimate $w(\xi,\tau) \leq C \langle
\tau - |\xi|^2 \rangle^{10}$.

It remains to establish the bilinear estimate (\ref{be}). Applying
\eqref{W-def} and monotonicity, we reduce to showing that
\begin{equation}\label{main-bilinear}
 \| \frac{w}{\langle \tau - |\xi|^2 \rangle} (\frac{f}{w} * \frac{g}{w}) \|_{Z^s}
\les \|f\|_{Z^s} \|g\|_{Z^s},
\end{equation}
for all non-negative smooth functions $f,g$.

Using \eqref{main-bilinear} together with the pasting Lemma, it
follows that \eqref{be} can be deduced by the bilinear estimates
stated in the following Proposition.
\begin{prop}\label{mainprop} Let $f,g$ be non-negative smooth functions and let $-1<s<0.$ Then,
\beq\label{YY-bilinear} \| \frac{w}{\langle \tau - |\xi|^2
\rangle} (\frac{f^{K^c}}{w} * \frac{g^{K^c}}{w}) \|_{Z^s} \lesssim
\|f^{K^c}\|_{Y^s} \|g^{K^c}\|_{Y^s},\eeq \beq \label{xxy}
||\frac{w}{\langle \tau-|\xi|^2 \rangle} (\frac{f^K}{w} *
\frac{g^K}{w})||_{Z^{s}} \lesssim ||f^K||_{\hat{X}^{s,\q,1}}
||g^K||_{\hat{X}^{s,\q,1}}, \eeq \beq \label{E3}
 \| \frac{\chi_{K} w}{\langle \tau - |\xi|^2 \rangle} (\frac{f^K}{w} * \frac{g^{K^c}}{w}) \|_{\hat{X}^{s,\q,1}}
\lesssim \|f^K\|_{\hat{X}^{s,\q,1}} \|g^{K^c}\|_{Y^{s}}, \eeq

\beq \label{E4}
 \| \frac{(1-\chi_{K}) w}{\langle \tau - |\xi|^2 \rangle} (\frac{f^K}{w} * \frac{g^{K^c}}{w}) \|_{Y^{s}}
\lesssim \|f^K\|_{\hat{X}^{s,\q,1}} \|g^{K^c}\|_{Y^{s}}. \eeq

\end{prop}

\begin{rem}\label{dropw}For $\tau=\tau_1 + \tau_2$, $\xi=\xi_1+\xi_2$ we have
$$ w(\xi,\tau) \leq C w(\xi_1,\tau_1) w(\xi_2,\tau_2);$$
hence we get the following pointwise estimate
\begin{equation}\label{k-point}
\frac{w}{\langle \tau - |\xi|^2 \rangle} (\frac{f}{w} *
\frac{g}{w}) \leq \frac{C}{\langle \tau - |\xi|^2 \rangle} (f*g),
\end{equation} which will turn out to be quite useful in the proof
of Proposition \ref{mainprop}.
\end{rem}

\section{The proof of the estimate \eqref{YY-bilinear}}

In this section we present the proof of the bilinear estimate
\eqref{YY-bilinear}. We start with a simple auxiliary lemma. Here
and henceforth, we will abbreviate $L^2_\xi L^2_\tau$ by $L^2$.
Analogously we denote by $L^\infty$ the mixed Lebesgue space
$L^\infty_\xi L^\infty_\tau.$

\begin{lem}Let $f,g$ be non-negative smooth functions. Then \begin{align}
 & ||\frac{1-\chi_{K}}{\langle \tau - |\xi|^2 \rangle}
(f*g)_{k,d}||_{L^{2}} \les 2^{-d} (2^{k+\frac{d}{2}})
||f||_{L^{2}} ||g||_{L^{2}}, \label{L22}\\
 &|| (f* \chi_{C_n}g)_{k,d}||_{
L^{2}_{\xi} L^{1}_{\tau}} \les  2^{k+\frac{n+d}{2}}
||f||_{L^{2}}||\chi_{C_n}g||_{L^{2}}\label{L21}.
\end{align}\end{lem}
\begin{proof} Let us prove \eqref{L22}. We have
\begin{align*}&||\frac{1-\chi_{K}}{\langle \tau - |\xi|^2 \rangle}
(f*g)_{k,d}||_{L^{2}} \leq \|(1-\chi_{K})\frac{\chi_{A_k \cap B_d
}}{\langle \tau - |\xi|^2
\rangle}\|_{L^{2}}\|f*g\|_{L^{\infty}}\\
&\leq  \|(1-\chi_{K})\frac{\chi_{A_k \cap B_d }}{\langle \tau -
|\xi|^2 \rangle}\|_{L^{2}}\|f\|_{L^{2}}\|g\|_{L^2}\end{align*}
where in the last step we have used Young's inequality. Moreover,
$$ \|(1-\chi_{K})\frac{\chi_{A_k \cap B_d }}{\langle \tau -
|\xi|^2 \rangle}\|_{L^{2}} \les 2^{-d} (2^{k+\frac{d}{2}})$$ since
$\langle \tau - |\xi|^2 \rangle \approx 2^{-d}$ on $B_d$ and the
variables $\xi$ and $\tau$ are confined to sets of measure
$2^{2k}$ and $2^d$ respectively.

We now prove \eqref{L21}. We have \begin{align}&|| (f*
\chi_{C_n}g)_{k,d}||_{ L^{2}_{\xi} L^{1}_{\tau}} \leq
\|\chi_{A_k\cap B_d}\|_{L^2}||f* \chi_{C_n}g||_{
L^{2}}\label{1}\\& \les 2^{k+\frac{d}{2}}||f||_{
L^{2}}||\chi_{C_n}g||_{ L_{\xi}^{2}L_{\tau}^1} \les
2^{k+\frac{d+n}{2}}||f||_{ L^{2}}||g||_{ L^{2}}\label{2}
\end{align}where to obtain \eqref{1} we used Cauchy-Schwartz in $\tau$,
while to obtain \eqref{2} we used that the variables $\xi$ and
$\tau$ are confined to sets of measure $2^{2k}$ and $2^d$
respectively together with Young's inequality. In the last step we
have then used Cauchy-Schwartz in $\tau$ and the fact that $\tau$
describes a set of measure approximately $2^{n}$.
\end{proof}

We are now ready to prove our bilinear estimate.

\begin{proof}[Proof of $\eqref{YY-bilinear}$] We wish to prove that

$$\| \frac{w}{\langle \tau - |\xi|^2 \rangle}
(\frac{f^{K^c}}{w} * \frac{g^{K^c}}{w}) \|_{Z^s} \lesssim
\|f^{K^c}\|_{Y^s} \|g^{K^c}\|_{Y^s}.$$

Define $h := f^{K^c}*g^{K^c}$. Then using the definition of the
space $Y^s$, the Pasting Lemma, and Remark \ref{dropw}, the
desired estimate will follow from the bilinear estimates:
\begin{equation}\label{yyx}
  \| \frac{\chi_K}{\langle \tau - |\xi|^2 \rangle}  h \|_{\hat{X}^{s,1/2,1}}
\lesssim \|f^{K^c}\|_{Y^s} \|g^{K^c}\|_{Y^s},
\end{equation}
\begin{equation}\label{yyl2}  \| \frac{1- \chi_K}{\langle \tau - |\xi|^2 \rangle} h \|_{\langle (\xi,\tau) \rangle^{s+1}L^2}
\lesssim \|f^{K^c}\|_{Y^s} \|g^{K^c}\|_{Y^s},
\end{equation}
\begin{equation}\label{yyl21}  \| \frac{w}{\langle \tau - |\xi|^2 \rangle} (1-\chi_K)
\left(\frac{f}{w} * \frac{g}{w}\right) \|_{\langle \xi
\rangle^{s}L^2_{\xi}L^1_\tau} \lesssim\|f^{K^c}\|_{Y^s}
\|g^{K^c}\|_{Y^s}.
\end{equation}
We divide our proof in three steps.

\

\noindent \textbf{Step 1}: Proof of \eqref{yyx}.

We decompose
$$h_k=\sum_{i,j}(f^{K^c}_i*g^{K^c}_j)_k,$$ where
$h_k=h_k(\xi,\tau),f_i=f_i(\xi_1,\tau_1)$ and
$g_j=g_j(\xi_2,\tau_2)$ with
$(\xi,\tau)=(\xi_1,\tau_1)+(\xi_2,\tau_2)$. Hence, in order for $
h_k$ to be non-zero, we must have that either of the following
cases holds:
\begin{itemize}
\item $j \leq i-11, |i-k| \leq 2$,
\item $i \leq j- 11, |j-k| \leq 2$,
\item $|i-j| \leq 10, k \leq \max\{i,j\} +1$.
\end{itemize}
Assume by symmetry that $i \leq j$. Then,
\begin{align}\label{sumfirst}
h_{k} = &\sum_{i \leq j-11;|j-k|\leq 2}
(f^{K^c}_{i}*g^{K^c}_{j})_{k} + \sum_{k\leq j+1; i\leq j\leq i+10} (f^{K^c}_{i}*g^{K^c}_{j})_{k}\\
\ & \nonumber \\
= &h'_k + h''_k\nonumber
\end{align}

A straightforward computation shows that,
\beq\label{1overtau}\|\frac{1}{\langle \tau - |\xi|^2
\rangle}\|_{L^2_{\xi}L^2_{\tau}(A_k \cap B_d)} \leq 2^{-d/2+k},
\eeq which will be used to analyze both $h'_k$ and $h''_k.$

\comment{We decompose $h'_k$ in the following way:
\begin{align*}h'_k &=  \sum_{j \leq k+5;|i-k|\leq 5;
|i-j|\leq 10} (f^{K^c}_{i}*g^{K^c}_{j})_{k}\\ &+\sum_{j \leq
k+5;|i-k|\leq 5; j < i-10} (f^{K^c}_{i}*g^{K^c}_{j})_{k} \\ &=
H_{11} + H_{12}.
\end{align*}
Let us estimate $H_{11}.$ From Young's inequality we have
\begin{align}\label{H11}\|H_{11}\|_{L^\infty_\tau L^\infty_\xi} & \leq \sum_{j \leq
k+5;|i-k|\leq 5;|i-j|\leq 10} 2^{-(i+j)(s+1)}\|f_i\|_{\langle \xi
\rangle^{s+1} L^2} \|g_j\|_{\langle \xi
\rangle^{s+1}L^2}\\
&\leq \sum_{j \leq k+5;|i-k|\leq 5; |i-j|\leq
10}2^{-(i+j)(s+1)}\|f_i\|_{Y^s}
\|g_j\|_{Y^s}\nonumber\end{align}where in the last inequality we
have used that $s+1>0.$}

\comment{Furthermore, using the identity:
$$\tau -|\xi|^2 = \tau_1 - |\xi_1|^2 + \tau_2 - |\xi_2|^2 -2 \xi_1
\cdot \xi_2,$$ together with the fact that $f^{K^c},g^{K^c}$ are
supported away from the parabola, and $|i-k|\leq 5, j \leq k+5$,
we get that $\chi_{B_d} H_{11}$ is supported on $d \geq 2k -10$.}

We start by estimating $h'_k$. By the definition of $\hat
X^{s,1/2,1}$ and of $K$, we get that $\displaystyle
\|\frac{\chi_{K}h'_k}{\langle \tau - |\xi|^2 \rangle}\|^{2}_{\hat
X^{s,1/2,1}}$ is bounded by

\begin{align*}&2^{2ks} \left( \sum_{d \leq 2k-4} \sum_{i \leq j-11;|j-k|\leq 2} 2^{\frac{d}{2}}\|\frac{1}{\langle \tau - |\xi|^2
\rangle}(f_i^{K^c} *g_j^{K^c})_{k,d}\|_{L^2}
\right)^2 \\
&\les \left( \sum_{d \leq 2k-4} \sum_{i \leq j-11;|j-k|\leq 2}
2^{\frac{d}{2}+ks}\sum_{d_1 \geq 2i-4; d_2 \geq
2j-4}\|\frac{1}{\langle \tau - |\xi|^2 \rangle}(f_{i,d_1}^{K^c}
*g_{j,d_2}^{K^c})_{k,d}\|_{L^2} \right)^2.\end{align*} In order to
estimate this last term, we notice that in light of the following
relation, $$ \tau-|\xi|^2=\tau_1 -|\xi_1|^2 + \tau_2 - |\xi_2|^2
-2 \xi_1 \cdot \xi_2, $$ non-trivial interactions happen only in
one of the following cases:
\begin{enumerate}
\item $|d - d_2|\leq 5, d_1 \leq d+6,$ \item  $|d -
d_1|\leq 5, d_2 \leq d+6,$ \item  $d_1,d_2 \geq d+7, |d_1 -
d_2|\leq 2.$
\end{enumerate}

Using \eqref{1overtau} together with H\"older's inequality, in all
such cases we have that the sum above is bounded by

$$\left( \sum_{d \leq 2k-4} \sum_{i \leq j-11;|j-k|\leq 2}
2^{\frac{d}{2}}2^{-\frac{d}{2}+k(s+1)}\sum_{d_1 \geq
2i-4}\sum_{d_2 \geq 2j-4}\|f_{i,d_1}^{K^c}
*g_{j,d_2}^{K^c}\|_{L^\infty} \right)^2$$

$$=  \left( \sum_{d \leq 2k-4} \sum_{i \leq j-11;|j-k|\leq 2}
2^{k(s+1)}\sum_{d_1 \geq 2i-4}\sum_{d_2 \geq
2j-4}\|f_{i,d_1}^{K^c} *g_{j,d_2}^{K^c}\|_{L^\infty} \right)^2 =
S_k.$$ In order to estimate $S_k$, in each of the cases
(i)-(ii)-(iii) we will use the following estimate, which follows
from  Young's inequality together with the definition of the space
$Y^s$ and the fact that $f^{K^c}$ and $g^{K^c}$ are supported away
from the parabola.
\begin{align}\label{H12}\|f_{i,d_1}^{K^c} *g_{j,d_2}^{K^c}\|_{L^\infty} & \leq
2^{-(d_1+d_2)\frac{(s+1)}{2}}\|f_{i,d_1}^{K^c}\|_{Y^s}
\|g_{j,d_2}^{K^c}\|_{Y^s}.\end{align}

\noindent \textit{Case} (i)-(ii). In this cases, the summation
over $d$ is finite, independent of $k$. Using \eqref{H12} and the
facts that $d_1 \geq 2i-4, d_2 \geq 2j-4,$ and $s+1>0$ we then get

\begin{align*}S_k \les & \left(\sum_{i \leq j-11;|j-k|\leq 2}
2^{(k-j-i)(s+1)}\|f^{K^c}_{i}\|_{Y^s}\|g^{K^c}_{j}\|_{Y^s}
\right)^2
\\ \les & \left(\sum_{i} 2^{-i(s+1)}\|f^{K^c}_{i}\|_{Y^s}\right)^2\sum_{|j-k|\leq
2}\|g^{K^c}_{j}\|^2_{Y^s}
\\ \les
&\|f^{K^c}\|^2_{Y^s}
\left(\sum_{p=-2}^{p=2}\|g^{K^c}_{p+k}\|^2_{Y^s}\right),\\
\end{align*}
where in the last inequality we used Cauchy-Schwartz inequality,
together with the fact that $s+1>0$. Notice that without this
assumption, we could not perform the desired summation.

Therefore, summing over all $k's$ we get:

 \beq\label{Sk(i)} \sum_{k}S_k \lesssim
\|f^{K^c}\|^2_{Y^s}\|g^{K^c}\|^2_{Y^s}.\eeq

\

\noindent \textit{Case} (iii). In this case the summation over $d$
is finite but dependent on $k$, while $|d_1-d_2|\leq 2$. Hence,
using the relations among the indexes together with the fact that
$\|f_{i,d_1}\|_{Y^s} \leq \|f^{K^c}_{i}\|_{Y^s},
\|g^{K^c}_{j,d_2}\|_{Y^s} \leq \|g_{j}\|_{Y^s}$, we get

\begin{align*}S_k \les & k^2\left(\sum_{i \leq j-11;|j-k|\leq 2} 2^{k(s+1)}\sum_{d_2 \geq 2j-4}
2^{-d_2(s+1)}\|f^{K^c}_{i}\|_{Y^s}\|g^{K^c}_{j}\|_{Y^s}\right)^2
\\\les & k^2\left(\sum_{i \leq j-11;|j-k|\leq 2}\sum_{d_2 \geq 2j-4}2^{(k-\frac{d_2}{2})(s+1)}
2^{-\frac{i}{2}(s+1)}2^{-\frac{k}{2}(s+1)}\|f^{K^c}_{i}\|_{Y^s}\|g^{K^c}_{j}\|_{Y^s}\right)^2\\
\les & k^2 2^{-k(s+1)}\|f^{K^c}\|_{Y^s}^2\|g^{K^c}\|_{Y^s}^2.
\end{align*}
Hence summing over all $k's$, under the assumption $s+1>0$ we get:

 \beq\label{Sk(iii)} \sum_{k}S_k \lesssim
\|f^{K^c}\|^2_{Y^s}\|g^{K^c}\|^2_{Y^s}.\eeq Combining
\eqref{Sk(i)} and \eqref{Sk(iii)}, we then obtain:

 \beq\label{h'k} \sum_{k} \|\frac{\chi_{K}h'_k}{\langle
\tau - |\xi|^2 \rangle}\|^{2}_{\hat X^{s,1/2,1}} \lesssim
\|f^{K^c}\|^2_{Y^s}\|g^{K^c}\|^2_{Y^s}.\eeq

 We now deal with
$h''_k$. From Young's inequality we have
\begin{align*}\|h''_k\|_{L^\infty} &\leq \sum_{k\leq j+1; i\leq j\leq i+10} 2^{-(i+j)(s+1)}\|f^{K^c}_i\|_{\langle \xi
\rangle^{s+1} L^2} \|g^{K^c}_j\|_{\langle \xi \rangle^{s+1}L^2} \\
&\leq \sum_{k\leq j+1; i\leq j\leq i+10}
2^{-(i+j)(s+1)}\|f^{K^c}_i\|_{Y^s}
\|g^{K^c}_j\|_{Y^s}.\end{align*}Hence, using \eqref{1overtau} and
the definition of $\hat X^{s,1/2,1}$ we get
\begin{align*}& \|\frac{\chi_{K}h''_k}{\langle \tau - |\xi|^2
\rangle}\|^{2}_{\hat X^{s,1/2,1}}  \\ & \leq 2^{2ks} \left(
\sum_{d \leq 2k-4} \sum_{k\leq j+1; i\leq j\leq i+10}
2^{\frac{d}{2}}2^{-\frac{d}{2}+k}2^{-(i+j)(s+1)}\|f^{K^c}_i\|_{Y^s}\|g^{K^c}_j\|_{Y^s}
\right)^2 \\ & \lesssim \left( \sum_{k\leq j+1; i\leq j\leq i+10}
k2^{(k-(i+j))(s+1)}\|f^{K^c}_i\|_{Y^s}\|g^{K^c}_j\|_{Y^s} \right)^2 \\
& \les \left(\sum_{k\leq j+1; i\leq j\leq i+10} k^2
2^{2(k-(i+j))(s+1)}\right)\|f^{K^c}\|^2_{Y^s}\|g^{K^c}\|^2_{Y^s},\end{align*}
where in the last inequality we have used Cauchy-Schwartz twice.
Since $s+1>0$, we have that \begin{align*} &\sum_k\sum_{k\leq j+1;
i\leq j\leq i+10} k^2 2^{2(k-(i+j))(s+1)} \les \sum_k\sum_{k\leq
j+1; i\leq j\leq i+10}
k^2 2^{-2i(s+1)}\\ & \les \sum_k  k^2 2^{-k(s+1)}\sum_{i }2^{-i(s+1)} \les 1.\\
\end{align*} Hence summing the inequality above over all $k$'s we
obtain \beq\label{h''k} \sum_{k}\|\frac{\chi_{K}h''_k}{\langle
\tau - |\xi|^2 \rangle}\|^{2}_{\hat X^{s,1/2,1}}  \lesssim
\|f^{K^c}\|^2_{Y^s}\|g^{K^c}\|^2_{Y^s}.\eeq Combining
\eqref{sumfirst},\eqref{h'k} and \eqref{h''k} we get the desired
estimate \eqref{yyx}.

We remark that similar computations show that,
 \beq\label{h''knew} \sum_{k}\|\frac{\chi_{K'}h_k}{\langle
\tau - |\xi|^2 \rangle}\|^{2}_{\hat X^{s,1/2,1}}  \lesssim
\|f^{K^c}\|^2_{Y^s}\|g^{K^c}\|^2_{Y^s},\eeq where
$$K'=\bigcup_{j} (A_j \cap B_{ \leq 2j+4}).$$

The reason is that for each $A_j$, $K'$ brings in few more modulations
than $K$, namely those between $2^{2j-3}$ and $2^{2j+4}$ and our argument 
easily tolerates this adjustment.

Hence, thanks to the Pasting Lemma, Lemma \ref{paste}, in the
proof of \eqref{yyl2} and \eqref{yyl21}, whenever $(\xi, \tau) \in
A_k \cap B_d$, we can restrict ourselves to the values $d \geq
2k+4$. A straightforward computation shows that in this regime we
have $2^{d-1} \leq |\tau| \leq 2^{d+2}$. Thus, in what follows, we
can assume \beq\label{tau2d}(\xi,\tau) \in A_k \cap B_d
\Rightarrow |\tau| \approx 2^d.\eeq

\

\noindent \textbf{Step 2.} Proof of  \eqref{yyl2}.

We use the following decomposition \beq (f^{K^c}*g^{K^c})_{k,d}
=\sum_{m,n} (\chi_{C_{m}} f^{K^c}* \chi_{C_{n}}
g^{K^c})_{k,d}.\eeq We observe that, under the assumption
\eqref{tau2d}, in order for the summands to be non-zero, one of
the following cases must hold:
\begin{itemize}
\item $n\leq d+5, |m-d| \leq 5$,
\item $m \leq d+5, |n-d| \leq 5$,
\item $m,n > d+5, |m-n| \leq 3$.
\end{itemize}

By symmetry we can then assume,
\begin{align}\label{sum3}
(f^{K^c}*g^{K^c})_{k,d}& =\sum_{n \leq d+5; |m-d|\leq 5}
(\chi_{C_{m}} f^{K^c}* \chi_{C_{n}} g^{K^c})_{k,d} +\\\nonumber &
+ \sum_{m,n \geq d+5;
|m-n|\leq 3}  (\chi_{C_{m}} f^{K^c}* \chi_{C_{n}} g^{K^c})_{k,d} \\
\nonumber & \\ \nonumber &= I_{k,d} + II_{k,d}.
\end{align}
From \eqref{L22} we get
$$||\frac{1-\chi_{K}}{\langle \tau - |\xi|^2 \rangle}
I_{k,d}||^2_{L^{2}} \leq 2^{-2d} (2^{2k+d})\sum_{n \leq d+5;
|m-d|\leq 5} ||\chi_{C_m}f^{K^c}||^2_{L^{2}}
||\chi_{C_n}g^{K^c}||^2_{L^{2}}.$$ Since on the complement of $K$,
$d \geq 2k -4$ we get (indeed we observed that we can restrict to
$d \geq 2k+4$),
\begin{align*}&||\frac{1-\chi_{K}}{\langle
\tau - |\xi|^2 \rangle} I_{k,d}||^2_{\langle (\xi,\tau)
\rangle^{s+1}L^{2}}\\ & \lesssim 2^{2k-d}\sum_{n \leq d+5;
|m-d|\leq 5}2^{(d-(m+n))(s+1)}
||\chi_{C_m}f^{K^c}||^2_{\la\tau\ra^{\frac{s+1}{2}}L^{2}}
||\chi_{C_n}g^{K^c}||^2_{\la\tau\ra^{\frac{s+1}{2}}L^{2}}\\
&\lesssim 2^{2k-d}\sum_{n \leq d+5; |m-d|\leq 5}2^{(d-(m+n))(s+1)}
||\chi_{C_m}f^{K^c}||^2_{Y^s}
||\chi_{C_n}g^{K^c}||^2_{Y^s},\end{align*} where in the last
inequality we used that $s+1>0.$
Therefore,
\begin{align*} & \sum_{d,k}||\frac{1-\chi_{K}}{\langle \tau - |\xi|^2
\rangle}I_{k,d}||^2_{\langle (\xi,\tau) \rangle^{s+1}L^{2}} \\
&\les \sum_{d} \left(\sum_{k \leq d/2 + 2} 2^{2k-d}\right)
\sum_{n\leq d+5;|m-d| \leq 5}2^{(d-(m+n))(s+1)}
||\chi_{C_m}f^{K^c}||^2_{Y^s} ||\chi_{C_n}g^{K^c}||^2_{Y^s}\\ &
\les ||f^{K^c}||^2_{Y^s} ||g^{K^c}||^2_{Y^s}.\end{align*} In order
to justify the last inequality, we can simplify the sum in $m$ by
taking $m=d$. Then we get \begin{align*}& \sum_d \sum_{n \leq d+5}
2^{-n(s+1)}
\|\chi_{C_d}f^{K^c}\|^2_{Y^s}\|\chi_{C_n}g^{K^c}\|^2_{Y^s} \\
&\leq \left(\sum_d \|\chi_{C_d}f^{K^c}\|^2_{Y^s}\right)
\left(\sum_{n} 2^{-n(s+1)}\right) \|g^{K^c}\|^2_{Y^s}\\& \les
\|f^{K^c}\|^2_{Y^s}\|g^{K^c}\|^2_{Y^s}\end{align*} where in the
last step we use the assumption $s+1>0$. We can therefore conclude
that
\begin{equation}\label{Ikd}
||\sum_{d,k}\frac{1-\chi_{K}}{\langle \tau - |\xi|^2 \rangle}
I_{k,d}||_{\langle (\xi,\tau) \rangle^{s+1}L^{2}}\les
\|f^{K^c}\|_{Y^s}\|g^{K^c}\|_{Y^s}.
\end{equation}
Analogously, using \eqref{L22} we get
\begin{align*}&\sum_{k,d}||\frac{1-\chi_{K}}{\langle \tau - |\xi|^2
\rangle}II_{k,d}||_{\langle (\xi,\tau) \rangle^{s+1}L^{2}} \\&
\les \sum_{d} \sum_{m,n> d+5;|m-n| \leq
3}2^{(d-(m+n))\frac{(s+1)}{2}} ||\chi_{C_m}f^{K^c}||_{Y^s}
||\chi_{C_n}g^{K^c}||_{Y^s} \\ &\leq \left(\sum_{d}\sum_{m,n >
d+5;|m-n|\leq
3}2^{(d-(m+n))\frac{(s+1)}{2}}\right)||f^{K^c}||_{Y^s}
||g^{K^c}||_{Y^s}\\& \les ||f^{K^c}||_{Y^s}
||g^{K^c}||_{Y^s}\end{align*} where in the last inequality we used
again that $s+1>0$. Combining the inequality above with
\eqref{Ikd} and \eqref{sum3} we obtain the desired estimate.

\

\noindent \textbf{Step 3.} Proof of \eqref{yyl21}.

Again, we use the decomposition \eqref{sumfirst}. We have the
following estimate
$$||f^{K^c}_{i}*g^{K^c}_{j}||_{L^{2}_{\xi}L^{1}_{\tau}}
\leq  ||f^{K^c}_i||_{L^{2}_{\xi}L^{1}_{\tau}}||g^{K^c}_j||_{L^{1}}
\leq 2^{j} ||f^{K^c}_{i}||_{L^{2}_{\xi}L^{1}_{\tau}}
||g^{K^c}_{j}||_{L^{2}_{\xi}L^{1}_{\tau}}.$$ Hence, since in the
support of $(1-\chi_{K})h_{k}$, $|\tau-|\xi|^{2}| \gtrsim 2^{2k}$,
we get \begin{align*}& ||\frac{1-\chi_{K}}{\langle \tau - |\xi|^2
\rangle}h'_{k} ||_{\langle \xi \rangle^{s} L^{2}_{\xi}
L^{1}_{\tau}} \\ &\leq \sum_{j \leq k+5; |k-i|\leq 5}2^{-2k} 2^{j}
2^{(k-j-i)s} ||f^{K^c}_{i}||_{\langle \xi \rangle^{s} L^{2}_{\xi}
L^{1}_{\tau}} ||g^{K^c}_{j}||_{\langle \xi \rangle^{s} L^{2}_{\xi}
L^{1}_{\tau}}.\end{align*} Therefore, square-summing in $k$ we get
\beq\label{h'kl1}\sum_k||\frac{1-\chi_{K}}{\langle \tau - |\xi|^2
\rangle}h'_{k} ||^2_{\langle \xi \rangle^{s} L^{2}_{\xi}
L^{1}_{\tau}} \les \|f^{K^c}\|^2_{Y^s}\|g^{K^c}\|^2_{Y^s}.\eeq In
order to justify such inequality, let us simplify the sum in $i$
by taking $i=k$. Then, using first Cauchy-Schwartz and then the
fact that $s+1>0$ we obtain \begin{align*}&\sum_k\left(\sum_{j
\leq k+5}2^{-2k+j(1-s)} ||f^{K^c}_{j}||_{\langle \xi \rangle^{s}
L^{2}_{\xi} L^{1}_{\tau}} ||g^{K^c}_{k}||_{\langle \xi \rangle^{s}
L^{2}_{\xi} L^{1}_{\tau}}\right)^2 \\ &\leq \sum_k \sum_{j\leq
k+5} 2^{-4(k-j)} ||f^{K^c}||^2_{Y^s} ||g^{K^c}_{k}||^2_{Y^s}\\
& \leq \sum_k \left(\sum_{p \geq -5} 2^{-4p}\right)
||f^{K^c}||^2_{Y^s} ||g^{K^c}_{k}||^2_{Y^s} \les
\|f^{K^c}\|^2_{Y^s}\|g^{K^c}\|^2_{Y^s}.\end{align*}In order to
control the behavior of $h''_k$, we perform a further
decomposition, that is \beq
 (f^{K^c}_{i}*g^{K^c}_{j})_{k,d} =\sum_{m,n} (\chi_{C_{m}} f^{K^c}_{i}* \chi_{C_{n}} g^{K^c}_{j})_{k,d}.\eeq
 Again, thanks to \eqref{tau2d}, in
order for the summands to be non-zero, one of the following cases
must hold:
\begin{itemize}
\item $m,n\leq d+5, |m-d| \leq 5$,
\item $m,n \leq d+5, |n-d| \leq 5$,
\item $m,n > d+5, |m-n| \leq 3$.
\end{itemize}
By symmetry we can then assume,
\begin{align}\label{sum2}
(f^{K^c}_{i}*g^{K^c}_{j})_{k,d}& =\sum_{m,n \leq d+5; |m-d|\leq 5}
(\chi_{C_{m}} f^{K^c}_{i}* \chi_{C_{n}} g^{K^c}_{j})_{k,d}
+\\\nonumber & + \sum_{m,n \geq d+5; |m-n|\leq 3} (\chi_{C_{m}}
f^{K^c}_{i}* \chi_{C_{n}} g^{K^c}_{j})_{k,d}= I_{i,j,k,d} +
II_{i,j,k,d}.
\end{align}
According to \eqref{2}, we have the following estimate
$$|| (\chi_{C_{m}} f^{K^c}_{i}* \chi_{C_{n}} g^{K^c}_{j})_{k,d}||_{
L^{2}_{\xi} L^{1}_{\tau}} \les  2^{k+\frac{n+d}{2}}
||\chi_{C_{m}}f^{K^c}_{i}||_{L^{2}}||\chi_{C_{n}}g^{K^c}_{j}||_{L^{2}}.$$
Hence, \begin{align*}&|| \frac{1-\chi_{K}}{\langle \tau - |\xi|^2
\rangle}(\chi_{C_{m}} f^{K^c}_{i}* \chi_{C_{n}}
g^{K^c}_{j})_{k,d}||_{ \langle \xi \rangle^{s} L^{2}_{\xi}
L^{1}_{\tau}} \\ &\les 2^{\frac{-d+n}{2}}2^{k(s+1)}
||\chi_{C_{m}}f^{K^c}_{i}||_{L^{2}}
||\chi_{C_{n}}g^{K^c}_{j}||_{L^{2}} \\
& \approx
2^{\frac{-d+n}{2}}2^{(k-(i+j))(s+1)}||\chi_{C_{m}}f^{K^c}_{i}||_{\langle
\xi \rangle^{s+1} L^{2}} ||\chi_{C_{n}}g^{K^c}_{j}||_{\langle \xi
\rangle^{s+1} L^{2}}\\ &\les
2^{\frac{-d+n}{2}}2^{(k-(i+j))(s+1)}||\chi_{C_{m}}f^{K^c}_{i}||_{Y^s}
||\chi_{C_{n}}g^{K^c}_{j}||_{Y^s}.\end{align*} Moreover,
$$\sum_{d} \sum_{m=d-5}^{d+5} \sum_{n \leq d+5} 2^{\frac{-d+n}{2}}
||\chi_{C_{m}}f^{K^c}_{i}||_{Y^s}
||\chi_{C_{n}}g^{K^c}_{j}||_{Y^s} \les ||f^{K^c}_{i}||_{Y^s}
||g^{K^c}_{j}||_{Y^s}.$$ To see this, we can simplify the sum in
$m$ and consider $m=d$. Then, the sum above becomes:
\begin{align*}&\sum_{d} \sum_{n \leq d+5} 2^{\frac{-d+n}{2}} ||\chi_
{C_{d}}f^{K^c}_{i}||_{Y^s}
||\chi_{C_{n}}g^{K^c}_{j}||_{Y^s}\\&=\sum_{p \geq -5}
2^{-\frac{p}{2}}  \sum_{d} ||\chi_{C_{d}}f^{K^c}_{i}||_{Y^s}
||\chi_{C_{p-d}}g^{K^c}_{j}||_{Y^s}.\end{align*} Applying
Cauchy-Schwartz first and then using the fact that $\sum_{p \geq
-5} 2^{-\frac{p}{2}}$ is bounded gives us the claim. Therefore, we
can conclude that
$$\sum_d ||\frac{1-\chi_{K}}{\langle \tau - |\xi|^2
\rangle}I_{i,j,k,d}||_{\langle \xi \rangle^{s} L^{2}_{\xi}
L^{1}_{\tau}} \les 2^{(k-(i+j))(s+1)}||f^{K^c}_{i}||_{Y^s}
||g^{K^c}_{j}||_{Y^s}.$$ Hence, \begin{align}\label{Iij} & \sum_k
\sum_{i,j\geq 5; |i-j|\leq 3}\sum_d ||\frac{1-\chi_{K}}{\langle
\tau - |\xi|^2 \rangle}I_{i,j,k,d}||_{\langle \xi \rangle^{s}
L^{2}_{\xi} L^{1}_{\tau}} \\\nonumber&\leq \sum_k \sum_{i,j\geq 5;
|i-j|\leq 3} 2^{(k-(i+j))(s+1)}||f^{K^c}_{i}||_{Y^s}
||g^{K^c}_{j}||_{Y^s} \les ||f^{K^c}||_{Y^s}
||g^{K^c}||_{Y^s}.\end{align} Now we treat the term involving
$II_{i,j,k,d}$. For this purpose we will need to use the weight
$w$. We have the estimate:
$$|| (\chi_{C_{m}} f^{K^c}_{i}* \chi_{C_{n}}
g^{K^c}_{j})_{k,d}||_{ L^{2}_{\xi} L^{1}_{\tau}} \leq  2^{k+d}
||\chi_{C_{m}}f^{K^c}_{i}||_{L^{2}}
||\chi_{C_{n}}g^{K^c}_{j}||_{L^{2}}.$$ If $(\xi_{1},\tau_{1}) \in
A_{i} \cap C_{m}$,  $(\xi_{2},\tau_{2}) \in A_{j} \cap C_{n}$ then
$(\xi_{1}+\xi_2,\tau_{1}+\tau_2) \in A_{k} \cap B_{d}$ only if
$\tau_{1}$ and $\tau_{2}$ have opposite sign. Therefore
\begin{align*}&|| \frac{w} {\langle \tau - |\xi|^2
\rangle}(1-\chi_{K}) (\frac{\chi_{C_{m}} f^{K^c}_{i}}{w}*
\frac{\chi_{C_{n}} g^{K^c}_{j}}{w})_{k,d}||_{ \langle \xi
\rangle^{s} L^{2}_{\xi} L^{1}_{\tau}} \\& \les 2^{(k-(i+j))(s+1)}
2^{10(d-m)} ||\frac{\chi_{C_{m}}f^{K^c}_{i}}{w}||_{\langle \xi
\rangle^{s+1}L^{2}} ||\frac{\chi_{C_{n}}g^{K^c}_{j}}{w}||_{\langle
\xi \rangle^{s+1}L^{2}}.\end{align*} As before, we can bound the
sum:
\begin{align*}& \sum_{d} \sum_{m,n \geq d+3;|m-n| \leq 3} 2^{10(d-m)}
||\frac{\chi_{C_{m}}f^{K^c}_{i}}{w}||_{\langle \xi
\rangle^{s+1}L^{2}} ||\frac{\chi_{C_{n}}g^{K^c}_{j}}{w}||_{\langle
\xi \rangle^{s+1}L^{2}}\\ &\les ||\frac{f^{K^c}_{i}}{w}||_{Y^s}
||\frac{g^{K^c}_{j}}{w}||_{Y^s}.\end{align*}
Hence,\begin{align}\label{IIij} &\sum_k \sum_{i,j\geq 5; |i-j|\leq
3}\sum_d ||\frac{1-\chi_{K}}{\langle \tau - |\xi|^2
\rangle}II_{i,j,k,d}||_{\langle \xi \rangle^{s} L^{2}_{\xi}
L^{1}_{\tau}}\\\nonumber & \les \sum_k \sum_{i,j\geq 5; |i-j|\leq
3} 2^{(k-(i+j))(s+1)}||\frac{f^{K^c}_{i}}{w}||_{Y^s}
||\frac{g^{K^c}_{j}}{w}||_{Y^s} \les ||f^{K^c}||_{Y^s}
||g^{K^c}||_{Y^s}.\end{align} Combining \eqref{Iij},\eqref{IIij},
with the decompositions \eqref{sumfirst},\eqref{sum2}, we obtain
the desired claim \eqref{yyl21}.

\end{proof}

\section{The proof of the estimate \eqref{xxy}}

This section is divided into two subsections. In the first
subsection we present a preliminary result from \cite{be}, and we
derive a new estimate similar to those in \cite{be}, which will be
used in the next subsection. There we exhibit the proof of the
bilinear estimate \eqref{xxy}.

\subsection{Interaction of Parabolas}

Here we describe how two ``parabolas" interact under convolution.
We will use these results in the proof of the estimate
\eqref{xxy}, where the two functions $f$ and $g$ are localized
near the parabola $\{\tau \approx |\xi|^2\}$ via the indicator
function $\chi_K$.

We remark that in the one-dimensional case, assuming that $(\xi,
\tau )= (\xi_1,\tau_1)+(\xi_2,\tau_2)$ one has the following
\textit{resonance estimate}
$$\max\{\la \tau - \xi^2 \ra,\la \tau_1 - \xi_1^2 \ra,\la \tau_2 - \xi_2^2 \ra\} \gtrsim \la\xi_1 \xi_2\ra.$$
Hence it follows immediately that if both the input frequencies
$\xi_1$ and $\xi_2$ are large, then it is not possible for all
three $(\xi,\tau), (\xi_1,\tau_1), (\xi_2,\tau_2)$ to lie close to
the parabola. This inequality is a powerful tool in the analysis
of the interaction of two functions $f,g$ supported near the
parabola. In higher dimensions the resonance inequality is no
longer true and we need to perform appropriate decompositions of
the frequency space, in order to investigate the interaction of
parabolas under convolution.

The substitute for the resonance inequality will be identity, \beq
\label{id} \tau-|\xi|^2=\tau_1 -|\xi_1|^2 + \tau_2 - |\xi_2|^2 -2
\xi_1 \cdot \xi_2. \eeq

We introduce now a few definitions. For each $c \in \R$ denote by
$$P_{c}= \{(\xi,\tau): \tau-|\xi|^{2}=c\}$$ and  by
$$\bar{P}_{c}=\{(\xi,\tau): \tau+|\xi|^{2}=c \}.$$ For notational
simplicity let  $P=P_{0}$ and $\bar{P}=\bar{P}_{0}$.

We denote by $\delta_{P_{c}}=\delta_{\tau-|\xi|^{2}=c}$ the
standard surface measure associated to the parabola $P_{c}$. Thus,
\beq \label{meas} \delta_{P_{c}}(f)=\int f d P_c=\int
f(\xi,|\xi|^2+c) \sqrt{1+4|\xi|^{2}} d\xi, \eeq for all smooth
functions $f$'s.

Then, the $L^2$ norm of the restriction of a smooth function $f$
to $P_{c}$ with respect to this measure, is given by

$$||f||_{L^{2}(P_{c})}=\left( \int f^{2}(\xi,|\xi|^{2}+c) \sqrt{1+4|\xi|^{2}} d\xi  \right)^{\q}$$

The following Lemma can be found in \cite{be} (see Propositions 1
and 2, section 4.1 there). Here we denote by  $P^{i}$ any of the
parabolas  $P_{c_{i}},\bar{P}_{c_{i}}, i=1,2$, for some constants
$c_1,c_2.$

\begin{lem} \label{pp1}
 Let $f \in L^{2}(P^1)$ and  $g \in L^{2}(P^2)$ be supported on  $A_i, A_j$ respectively. Then, \beq \label{ge} ||f \delta_{P^{1}} *
g \delta_{P^{2}} ||_{L^{2}} \leq 2^{\min{(i,j)}}
||f||_{L^{2}(P^{1})}||g||_{L^{2}(P^{2})}. \eeq

 Moreover, if $i \leq j$, $|c_{1}| \leq 2^{2i-2}$ and
$|c_{2}| \leq 2^{2j-2}$,
 then \beq \label{ge2} ||f \delta_{P^{1}} * g \delta_{P^{2}}
||_{L^{2}(|(\xi,\tau)| \approx 2^{j},|\tau-|\xi|^{2}| \leq d)}
\lesssim d^{\q} ||f||_{L^{2}(P^{1})}||g||_{L^{2}(P^{2})}. \eeq

\end{lem}

The next result is needed due to our particular choice of spaces.

\begin{lem}
 Let $f \in L^{2}(P^1)$ and  $g \in L^{2}(P^2)$ be supported on  $A_i, A_j$ respectively. Assume that
$|c_{1}| \leq 2^{2i+10}$ and $|c_{2}|\leq 2^{2j+10}$. Then, \beq
\label{c2} ||f \delta_{P^{1}} * g \delta_{P^{2}}
||_{L^{2}_{\xi}L^{1}_{\tau}(|\xi| \approx 2^{k})} \lesssim
2^{k+\frac{i+j}{2}} ||f||_{L^{2}(P^{1})}||g||_{L^{2}(P^{2})}. \eeq

\end{lem}

\begin{proof}Without loss of generality we can assume $c_{1}=c_{2}=0$.
Also, let us consider the case when $P^1=P$ and $P^2=P$. The other
cases can be treated in a similar way.

Assuming that $f \delta_{P} * g \delta_{P} \in
L^{2}_{\xi}L^{1}_{\tau}(|\xi| \approx 2^{k})$, then its norm is
controlled by estimating $|(f \delta_{P} * g \delta_{P})h|$ for
any $h \in L^{2}_{\xi}L^{\infty}_{\tau}$ supported at frequency
$2^{k}$. For any such $h$ we have:

$$(f \delta_{P}* g \delta_{P})h=
\int f(\xi) g(\eta) h(\xi+\eta, |\xi|^2+|\eta|^2)
\sqrt{1+4|\xi|^2}  \sqrt{1+4|\eta|^2} d\xi d\eta. $$ We decompose
$\R^{2} = \cup_{l \in \Z^{2}} Q_{l}^{k}$ in cubes $Q^{k}_{l}$
centered at $2^{k}l$ and of size $2^{k}$. Then we split

$$f=\sum_{l \in \Z^{2}} f_{l}  \ \ \ \mbox{and} \ \ \ g=\sum_{l \in \Z^{2}} g_{l}$$
where $f_{l}$ is the part of $f$ localized in $Q_{l}^{k}$ and
similarly for $g_l$. Since $h$ is supported at frequency $2^{k}$
we obtain,

$$(f \delta_{P} * g \delta_{P})h=\sum_{l \in \Z^{2}} \int f_{l}(\xi) g_{-l}(\eta) h(\xi+\eta,
|\xi|^2+|\eta|^2) \sqrt{1+4|\xi|^2}  \sqrt{1+4|\eta|^2} d\xi
d\eta.
$$
For fixed $l$ we evaluate:

\begin{align*}|(f_{l} \delta_{P} * g_{-l} \delta_{P})h| &\leq ||f_{l}||_{L^{2}(P)}
||g_{-l}||_{L^{2}(P)}\\
 &\times \left( \int h^{2}(\xi+\eta,|\xi|^2+|\eta|^2) \sqrt{1+4|\xi|^2}  \sqrt{1+4|\eta|^2} d\xi d\eta
 \right)^{\q}.\end{align*}
Since $h \in L^{2}_{\xi}L^{\infty}_{\tau}$ we can suppose that we
have $\tilde{h} \in L^{2}_{\xi}$ and estimate:

$$\int \tilde{h}^{2}(\xi+\eta) \sqrt{1+4|\xi|^2}  \sqrt{1+4|\eta|^2} d\xi d\eta \leq 2^{2k+i+j} ||\tilde{h}||^{2}_{L^{2}}$$
from which we can conclude
$$|(f_{l} \delta_{P}* g_{-l}\delta_{P})h| \leq 2^{k+\frac{i+j}{2}} ||f_{l}||_{L^{2}(P)} ||g_{-l}||_{L^{2}(P)}
||h||_{L^{2}_{\xi}L^{\infty}_{\tau}}.$$ Summing up with respect to
$l$ gives us the claim in (\ref{c2}), as long as  $f \delta_{P} *
g \delta_{P} \in L^{2}_{\xi}L^{1}_{\tau}(|\xi| \approx 2^{k})$. In
order to guarantee this last fact, we can pick $h \in
L^{2}_{\xi,\tau}$ and perform a similar computation to the one
above. In this way we obtain an $L^{2}$ estimate for $f
\delta_{P}* g \delta_{P}$. However, $ L^{2}_{\xi,\tau}(|\xi|
\approx 2^{k}) \subset L^{2}_{\xi}L^{1}_{\tau}(|\xi| \approx
2^{k})$, and this concludes the proof.
\end{proof}

\subsection{Proof of the bilinear estimate \eqref{xxy}.} We start with an auxiliary lemma.

\begin{lem} \label{u1}Let $f,g$ be non-negative smooth functions. Assume that $0 \leq i
\leq j, d_1 \leq 2i-4$, $d_2 \leq 2j-4$. Then, for $|k-j| \leq 5$ and $d_3 \leq 2k-4$, the following estimates
hold.

\beq \label{b2}
||(f_{i,d_1}*g_{j,d_2})_{k,d_3}||_{\hat{X}^{0,-\q}} \les
2^{-\frac{i+j}{2}}  ||f_{i,d_1}||_{\hat{X}^{0,\q}}
||g_{j,d_2}||_{\hat{X}^{0,\q}}. \eeq

Moreover, if $|i-j| \leq 2$, $k \leq j+2$, $d_1 \leq 2i-4$ and $d_2 \leq 2j-4$, then

\beq \label{b10} ||f_{i,d_1}*g_{j,d_{2}}||_{L^{2}} \leq
||f_{i,d_1}||_{\hat{X}^{0,\q}} ||g_{j,d_{2}}||_{\hat{X}^{0,\q}},
\eeq

\beq \label{b11}
||(f_{i,d_1}*g_{j,d_{2}})_{k}||_{L^{2}_{\xi}L^{1}_{\tau}} \les
2^{k} ||f_{i,d_1}||_{\hat{X}^{0,\q}}
||g_{j,d_{2}}||_{\hat{X}^{0,\q}}. \eeq

\end{lem}

\begin{proof}

The proof of \eqref{b2} and \eqref{b10} can be found in
\cite{be}, see the proof of Proposition 3, section 4.2. These
estimates can be also derived directly from the statements of
Proposition 3, section 4.2 in \cite{be}.

Thus, we are left with the proof of \eqref{b11}, which is specific
to our problem. We aim to apply the results in the Subsection 4.1.
Notice that on each $A_i$ the parabolas $P_{c}$ make an angle of
approximately $2^{-i}$ with the $\tau$ axis. Thus, recalling
\eqref{meas}, we have the following relation between measures:

$$d\xi d\tau \approx 2^{-i} dP_{c} dc.$$
Therefore for each $h$ supported in $A_{k} \cap B_d$ we have \beq
\label{e100} ||h||^{2}_{\hat{X}^{0,\pm \q}} \approx 2^{\pm d}
\int_{2^{d-1}}^{2^{d+1}} ||h||^{2}_{L^{2}(P_{b})} 2^{-k} db. \eeq
For notational simplicity, we denote $f=f_{i,d_1}, g=g_{j,d_2}$.
We apply (\ref{c2}) and \eqref{e100} to evaluate the following
 norm,

\begin{align*}&||f * g||_{L^{2}_{\xi} L^{1}_{\tau}} \leq \int_{I_{1}} \int_{I_{2}} ||f
\delta_{P_{b_{1}}} * g \delta_{P_{b_{2}}} ||_{L^{2}_{\xi} L^{1}_{\tau}}  2^{-i-j}
db_{1} db_{2} \leq \\
&\int_{I_{1}} \int_{I_{2}} 2^{k-\frac{i+j}{2}} ||f||_{L^{2}(P_{b_{1}})}
||g||_{L^{2}(P_{b_{2}})} db_{1} db_{2} \leq \\
&2^{k} \left( \int_{I_{1}} (1+b_{1})^{-1} db_{1}
\right)^{\q} ||f||_{\hat{X}^{0,\q}}
 \left( \int_{I_{2}} (1+b_{2} )^{-1} db_{2} \right)^{\q} ||g||_{\hat{X}^{0,\q}} \approx
 \\
&2^{k} ||f||_{\hat{X}^{0,\q}}
||g||_{\hat{X}^{0,\q}}.\end{align*} Here we used the
fact that $I_{1} \approx [2^{d_{1}-1},2^{d_{1}+1}],$ which gives
us $\int (1+b_{1})^{-1} db_{1} \approx 1 $. Same thing for the
integral with respect to $b_{2}$.

We remark that the proof of \eqref{b2}-\eqref{b10} follows the
same lines as the proof above, and uses Lemma \ref{pp1}, which we
have thus stated for convenience of the reader. However, the proof
is slightly more involved, hence we prefer to refer the reader to
\cite{be} as well.

\end{proof}

\begin{proof}[Proof of $\eqref{xxy}$] We wish to prove that
$$||\frac{w}{\langle \tau-|\xi|^2 \rangle} (\frac{f^K}{w} *
\frac{g^K}{w})||_{Z^{s}} \lesssim ||f^K||_{\hat{X}^{s,\q,1}}
||g^K||_{\hat{X}^{s,\q,1}}.
$$

On behalf of (\ref{k-point}) we reduce to showing the following
bound, \beq\label{dani}||\frac{1}{\langle \tau-|\xi|^2 \rangle}
(f^K * g^K)||_{Z^{s}} \lesssim ||f^K||_{\hat{X}^{s,\q,1}}
||g^K||_{\hat{X}^{s,\q,1}}. \eeq

Now, we assume by symmetry that $i \leq j$; we perform the
following decomposition \beq \label{dec11} (f^K_{i}*
g^K_{j})=\chi_{K}(f^K_{i}* g^K_{j})+(1-\chi_{K})(f^K_{i}*
g^K_{j}). \eeq We continue with decomposing even further the first
term, that is \beq \label{dec111}\chi_{K}(f^K_{i}* g^K_{j}) =
\sum_{k} \sum_{d_3 \leq 2k-4} \sum_{d_{1} \leq 2i-4}
 \sum_{d_2 \leq 2j-4} (f_{i,d_1}*g_{j,d_2})_{k,d_3}.\eeq
By the definition of $K$, in the support of $f^K_i$ we have $\tau
\geq 2^{2i-4}$ while in the support of $g^K_{j}$ we have $\tau
\geq 2^{2j-4}$. Hence, in the support of $f^K_{i}*g^K_j$ we have
$\tau \geq 2^{2j-4}$ and in the support of
$\chi_{K}(f^K_{i}*g^K_j)$ we have $|\xi| \geq 2^{j-5}$. This means
we have  nontrivial interactions only in the case $|k-j| \leq 5$.

Hence, using (\ref{b2}) together with \eqref{dec111}, we get the
following estimate

\begin{align*} &||\frac{\chi_{K}}{\langle \tau-|\xi|^2 \rangle} (f^K_{i} *
g^K_{j})||_{\hat{X}^{0,\q,1}} \approx ||\chi_{K} (f^K_{i} *
g^K_{j})||_{\hat{X}^{0,-\q,1}}\\ & \les
 \sum_{|k-j| \leq 5} \sum_{d_3 \leq 2k-4} \sum_{d_{1} \leq 2i-4}
\sum_{d_2 \leq 2j-4} | |(f_{i,d_{1}} *
g_{j,d_2})_{k,d_{3}}||_{\hat{X}^{0,-\q,1}} \\ & \les \sum_{|k-j|
\leq 5} \sum_{d_3 \leq 2k-4} \sum_{d_{1} \leq 2i-4} \sum_{d_2 \leq
2j-4} 2^{-\frac{i+j}{2}} ||f_{i,d_{1}}||_{\hat{X}^{0,\q}}
||g_{j,d_{2}}||_{\hat{X}^{0,\q}}
\\ & \les 2^{-\frac{i+j}{2}} \sum_{|k-j| \leq 5} \sum_{d_3 \leq
2k-4}
||f^K_{i}||_{\hat{X}^{0,\q,1}} ||g^K_{j}||_{\hat{X}^{0,\q,1}}\\
& \les j 2^{-\frac{i+j}{2}} ||f^K_{i}||_{\hat{X}^{0,\q,1}}
||g^K_{j}||_{\hat{X}^{0,\q,1}}.\end{align*} For general $s$ this
becomes
$$||\frac{\chi_{K}}{\langle \tau-|\xi|^2 \rangle} (f^K_{i} * g^K_{j})||_{\hat{X}^{s,\q,1}}
\leq j 2^{-(1+s)i} 2^{\frac{i-j}{2}}
||f^K_{i}||_{\hat{X}^{s,\q,1}} ||g^K_{j}||_{\hat{X}^{s,\q,1}}.$$
Summing up with  respect to $i,j$ gives us

\begin{align} \label{k1} & ||\frac{\chi_{K}}{\langle \tau-|\xi|^2 \rangle}
(f^K * g^K)||_{X^{s,\q,1}}  \\ & \les \sum_{i \leq j} j
2^{-(1+s)i} 2^{\frac{i-j}{2}} ||f^K_{i}||_{\hat{X}^{s,\q,1}}
||g^K_{j}||_{\hat{X}^{s,\q,1}} + \mbox{symmetric term} \nonumber
\\ & \lesssim ||f^K||_{\hat{X}^{s,\q,1}}
||g^K||_{\hat{X}^{s,\q,1}}. \nonumber \end{align} Finally, since
$\hat{X}^{s,\q,1}$ and $Z^s$ paste nicely in the set $K$, we have
obtained that \beq\label{xxy1}||\frac{\chi_{K}}{\langle
\tau-|\xi|^2 \rangle} (f^K* g^K)||_{Z^s} \lesssim
||f^K||_{\hat{X}^{s,\q,1}} ||g^K||_{\hat{X}^{s,\q,1}}.\eeq

For the second term in (\ref{dec11}) we make the following
observation. If $i \leq j-5$, then $f^K_{i} * g^K_{j}$ is
supported in $\cup_{k=j-1}^{j+1} A_{k} \cap B_{\leq 2j-3}$. Hence,
by the pasting Lemma, in this case the $Z^s$ norm of
$(1-\chi_{K})(f^K_{i}
* g^K_{j})$ is controlled by the $\hat{X}^{s,\q,1}$ norm. Morally, in this case we keep
the interaction close enough to $P$ so that we can treat the
estimate in the same way as the previous one. Indeed, using the
same computations as above we obtain, \beq \label{k2} \sum_{|i-j|
\geq 5}||\frac{1-\chi_{K}}{\langle \tau - |\xi|^2 \rangle}
(f^K_{i} * g^K_{j})||_{\hat{X}^{s,\q,1}} \les
||f^K||_{\hat{X}^{s,\q,1}} ||g^K||_{\hat{X}^{s,\q,1}}.
 \eeq In the
case when $|i-j| \leq 4$ we perform the following decomposition,
\beq \label{dec20}
 (1-\chi_{K})(f^K_{i} * g^K_{j}) = \sum_{k} \sum_{d_{1} \leq 2i-4}
  \sum_{d_2 \leq 2j-4} (1-\chi_{K}) (f_{i,d_1} * g_{j,d_2})_{k}.
\eeq We notice that $f^K_{i} * g^K_{j}$ is supported in a region
where
  $|\xi| \leq 2^{j+10}$, $\tau \geq 2^{2j-4}$ and $|\tau-|\xi|^2| \geq 2^{2j-10}
  $.
Using the decomposition \eqref{dec20} and the estimate (\ref{b10})
we obtain,
\begin{align} \label{k3} & \sum_{|i-j| \leq 4}||\frac{1-\chi_{K}}{\langle \tau-|\xi|^2 \rangle}
(f^K_i* g^K_j)||_{\langle (\xi,\tau)\rangle^{s+1} L^{2}} \\ & \les
\sum_{|i-j| \leq 4} \sum_{k \leq j+10} 2^{(s+1)(k-2j)}
||f^K_{i}||_{\hat{X}^{s,\q,1}} ||g^K_{j}||_{\hat{X}^{s,\q,1}}
\nonumber
\\ & \lesssim ||f^K||_{\hat{X}^{s,\q,1}}
||g^K||_{\hat{X}^{s,\q,1}},\nonumber \end{align} while using the
estimate (\ref{b11}) we get,
\begin{align} \label{k4} & \sum_{|i-j| \leq 4}||\frac{1-\chi_{K}}{\langle \tau-|\xi|^2 \rangle}
(f^K_i * g^K_j)||_{\langle \xi \rangle^{s}
L^{2}_{\xi}L^{1}_{\tau}}
\\&\les \sum_{|i-j| \leq 4} \sum_{k \leq j+10} 2^{(s+1)(k-2j)} ||f^K_{i}||_{\hat{X}^{s,\q,1}} ||g^K_{j}||_{\hat{X}^{s,\q,1}} \nonumber \\
& \lesssim||f^K||_{\hat{X}^{s,\q,1}}
||g^K||_{\hat{X}^{s,\q,1}}.\nonumber\end{align}

The estimates \eqref{k2},\eqref{k3} and \eqref{k4}, the definition
of the space $Y^s$, and the fact that $Y^s$ and $Z^s$ paste nicely
outside of $K$ then imply

\beq\label{xxy2}||\frac{1-\chi_{K}}{\langle \tau-|\xi|^2 \rangle}
(f^K
* g^K)||_{Z^s} \lesssim
||f^K||_{\hat{X}^{s,\q,1}} ||g^K||_{\hat{X}^{s,\q,1}}.\eeq
Combining \eqref{xxy1}  and \eqref{xxy2} we get the claim in
(\ref{dani}).

\end{proof}

\section{The proof of the estimate \eqref{E3}}

In this section we present the proof of the bilinear estimate
\eqref{E3}. We start by proving an auxiliary lemma.

\begin{lem}Let $f,g$ be non-negative smooth functions. If  $k \geq i-10$, $d_{3} \leq
2k-4$ and $d_1 \leq 2i-4$ then,

\beq \label{b20} ||(f_{i,d_1}*g)_{k,d_3}||_{\hat{X}^{0,-\q}} \lesssim
2^{\frac{i-k}{2}}||f_{i,d_1}||_{\hat{X}^{0,\q}} ||g||_{L^{2}}.
\eeq

If $d_{2} \leq 2i-10$, $d_{3} \leq 2k-4$ and $d_1 \leq 2i-4$, then

\beq \label{b21} ||(f_{i,d_1}*g_{d_2}^{K^c})_{k,d_3}||_{\hat{X}^{0,-\q}} \lesssim
2^{\frac{3d_2-2i-2d_3}{4}}||f_{i,d_1}||_{\hat{X}^{0,\q}}
||g^{K^c}_{d_2}||_{L^{2}}. \eeq

If $k \leq i-10$ and $d_3 \leq 2k-4$ and $d_1 \leq 2i-4$, then

\beq \label{b22} ||(f_{i,d_1} * \frac{g}{w})_{k,d_3}||_{\hat{X}^{0,-\q}}
\lesssim 2^{-4i} ||f_{i,d_1}||_{\hat{X}^{0,\q}} ||g||_{L^{2}}.
\eeq

\end{lem}

\begin{proof}

By duality, (\ref{b20}) is equivalent to:

$$||f_{i,d_1} * g_{k,d_3}||_{L^{2}} \lesssim 2^{\frac{i-k}{2}} ||f_{i,d_1}||_{\hat{X}^{0,\q}} ||g_{k,d_3}||_{\hat{\bar{X}}^{0,\q}}$$
which can be obtain by a similar argument to the one for
(\ref{b2}), in light of the fact that (\ref{ge}) allows us to work
with either $\delta_{P_{c}}$ or $\delta_{\bar{P}_{c}}$. For a
detailed argument we refer to \cite{be}.

 For the proof of (\ref{b21}) we split $A_i \cap
B_{d_1}=\cup_{\a} D_{\a}$, where $D_{\a}$ are disjoint sets of
sizes $2^{-i}2^{d_1} \times 2^{\frac{d_2}{2}} \times 2^{d_2}$, the
first size being in the direction of $n_{\a}$ and the last one in
the direction of $\tau$. Here by $n_{\a}$ we mean  one of the
normal directions to $P_{d_2} \cap D_{\a}$.

We replace $f_{i,d_1}$ by $f$, so that we do not carry all the indexes.

Then, if $f_{\a}$ is the part of $f$ localized in $D_{\a}$, the
$f_{\a} * g $'s have, essentially, disjoint support with respect
to $\a$. Hence,

$$||f*g||_{L^{2}}^{2} \approx \sum_{\a} ||f_{\a} * g||_{L^{2}}^{2} \leq \sum_{\a} ||f_{\a}||_{L^{1}}^{2} ||g||^{2}_{L^{2}} \lesssim $$

$$2^{-i+d_1} 2^{\frac{3}{2}d_2} \sum_{\a} ||f_{\a}||_{L^{2}}^{2} ||g||^{2}_{L^{2}}
\lesssim 2^{-i+d_1} 2^{\frac{3}{2}d_2}||f||_{L^{2}}^{2}
||g||^{2}_{L^{2}}.$$

Thus,
\begin{align*}||(f * g)_{j,d_3}||_{\hat{X}^{0,-\q}} & \lesssim 2^{-\frac{d_3}{2}} ||f*g||_{L^{2}}
 \lesssim\ddot{}  2^{\frac{3d_2-2i-2d_3}{4}} 2^{\frac{d_1}{2}} ||f||_{L^{2}} ||g||_{L^{2}}
 \\& \approx 2^{\frac{3d_2-2i-2d_3}{4}} ||f_{i,d_1}||_{\hat{X}^{0,\q}} ||g||_{L^{2}},\end{align*}
that is the desired estimate.

 We start the proof of (\ref{b22}) with a
geometrical observation. The function $f$ is supported in a region
where $\tau \geq 2^{2i-4}$ and we want to restrict $f * g$ in a
region where $\tau \leq 2^{2k+4}$. This can be achieved only by
restricting the support of $g$ in a region where $\tau \leq
-2^{2i-5}$. This allows us to obtain pretty loose estimates,
thanks to the weight $w$.

Indeed, we can run the $L^{1}
* L^{2} \ri L^{2}$ argument to get
$$||(f * \frac{g}{w})_{k,d_{3}}||_{\hat{X}^{0,-\q}} \les ||f * \frac{g}{w}||_{L^{2}}
\les ||f||_{L^{1}} ||\frac{g}{w}||_{L^{2}} \les 2^{i}
||f_{i,d_1}||_{\hat{X}^{0,\q}} 2^{-20i}||g||_{L^{2}}.$$ This
concludes our proof.
\end{proof}

\begin{proof}[Proof of $(\ref{E3})$] We wish to prove that,
\beq \label{E3new}
 \| \frac{\chi_{K} w}{\langle \tau - |\xi|^2 \rangle} (\frac{f^K}{w} * \frac{g^{K^c}}{w}) \|_{\hat{X}^{s,\q,1}}
\lesssim \|f^K\|_{\hat{X}^{s,\q,1}} \|g^{K^c}\|_{Y^{s}}. \eeq

Since $w=1$ in $K$, we can drop $w$ from $f^K/w$. Then we
decompose
\begin{align*}&\chi_{K} (f^K *\frac{g^{K^c}}{w})= \sum_{k} \sum_{d_3 \leq
2k-4} \sum_{i} \sum_{d_1 \leq 2i-4} \sum_{d_2}  (f_{i,d_1}
* \frac{g^{K^c}_{d_2}}{w})_{k,d_3}\\ &=
\sum_{k} \sum_{d_3 \leq 2k-4} \sum_{i \leq k+10} \sum_{d_1 \leq
2i-4} \sum_{d_2} (f_{i,d_1}
* \frac{g^{K^c}_{d_2}}{w})_{k,d_3}\\ &+ \sum_{k} \sum_{d_3 \leq 2k-4} \sum_{i \geq
k+11} \sum_{d_1 \leq 2i-4} \sum_{d_2}(f_{i,d_1}
* \frac{g^{K^c}_{d_2}}{w})_{k,d_3}\\
&=I_{1}+I_{2}\end{align*} We decompose $I_{1}$ even further
\begin{align*}I_{1}&=\sum_{k} \sum_{d_3 \leq 2k-4} \sum_{i \leq k+10}
\sum_{d_1 \leq 2i-4} \sum_{d_2 \geq \frac{i}{10}} (f_{i,d_1}
* \frac{g^{K^c}_{d_2}}{w})_{k,d_3}\\&+\sum_{k}
\sum_{d_3 \leq 2k-4} \sum_{i \leq k+10} \sum_{d_1 \leq 2i-4}
\sum_{d_2 < \frac{i}{10}}(f_{i,d_1}
* \frac{g^{K^c}_{d_2}}{w})_{k,d_3}\\ &=I_{11}+I_{12}\end{align*}

We use (\ref{b20}) to estimate $I_{11}$:

\begin{align*}&||\frac{1}{\la \tau - |\xi|^2 \ra}I_{11}||_{X^{s,\q,1}} \les ||I_{11}||_{X^{s,-\q,1}} \\&\les
\sum_{k} \sum_{d_3 \leq 2k-4} \sum_{i \leq k+10} \sum_{d_1 \leq
2i-4} \sum_{d_2 \geq \frac{i}{10}} 2^{ks} ||(f_{i,d_1} *
g_{d_2}^{K^c})_{k,d_3}||_{\hat{X}^{0,-\q}} \\ &\les \sum_{k}
\sum_{d_3 \leq 2k-4} \sum_{i \leq k+10} \sum_{d_1 \leq 2i-4}
\sum_{d_2 \geq \frac{i}{10}} 2^{ks} 2^{\frac{i-k}{2}}
||f_{i,d_1}||_{\hat{X}^{0,\q}} ||g_{d_2}^{K^c}||_{L^{2}}
\\&\les
\sum_{k} \sum_{d_3 \leq 2k-4} \sum_{i \leq k+10} \sum_{d_2 \geq
\frac{i}{10}} 2^{ks} 2^{\frac{i-k}{2}} 2^{-\frac{(s+1)d_2}{2}}
||f_{i}^K||_{\hat{X}^{0,\q,1}} ||g_{d_2}^{K^c}||_{Y^{s}} \\
&\les\sum_{k} \sum_{d_3 \leq 2k-4} \sum_{i \leq k+10} 2^{(k-i)s}
2^{\frac{i-k}{2}} 2^{-(s+1)\frac{i}{20}}
||f_{i}^K||_{\hat{X}^{s,\q,1}} ||\chi_{B_{\geq \frac{i}{10}}}
g^{K^c}||_{Y^{s}} \\ &\les \sum_{k} \sum_{i \leq k+10}  k
2^{(k-i)(s-\q)} 2^{-(s+1)\frac{i}{20}}
||f_{i}^K||_{\hat{X}^{s,\q,1}} ||g^{K^c}||_{Y^{s}} \les
||f^K||_{\hat{X}^{s,\q,1}} ||g^{K^c}||_{Y^{s}}.\end{align*} When
estimating $I_{12}$ we notice that unless $|k-i| \leq 2$ we have
trivial estimates. Indeed the support of $\chi_{B_{\leq
\frac{i}{10}}} g^{K^c}$ cannot move, via convolution, the support
of $f_i$ too much. We use (\ref{b21}) and compute:
\begin{align*}&||\frac{1}{\la \tau - |\xi|^2 \ra}I_{12}||_{X^{s,\q,1}} \les ||I_{12}||_{X^{s,-\q,1}}\\& \les
\sum_{k} \sum_{d_3 \leq 2k-4} \sum_{i = k-2}^{k+2} \sum_{d_1 \leq
2i-4} \sum_{d_2 \leq \frac{i}{10}} 2^{ks} ||(f_{i,d_1} *
g_{d_2}^{K^c})_{k,d_3}||_{\hat{X}^{0,-\q}} \\ &\les \sum_{k}
\sum_{d_3 \leq 2k-4} \sum_{i = k-2}^{k+2} \sum_{d_1 \leq 2i-4}
\sum_{d_2 \leq \frac{i}{10}} 2^{ks}
2^{\frac{3d_2-2i-2d_3}{4}}||f_{i,d_1}||_{\hat{X}^{0,\q}}
||g_{d_2}^{K^c}||_{L^{2}} \\ & \les \sum_{k} \sum_{d_3 \leq 2k-4}
\sum_{i = k-2}^{k+2} \sum_{d_2 \leq \frac{i}{10}}
2^{\frac{3d_2-2i-2d_3}{4}} 2^{-\frac{(s+1)d_2}{2}}
||f^K_{i}||_{\hat{X}^{s,\q,1}} ||g_{d_2}^{K^c}||_{Y^{s}} \\ & \les
\sum_{k} \sum_{d_3 \leq 2k-4} \sum_{i = k-2}^{k+2}
2^{\frac{-i-2d_3}{4}} ||f_{i}^K||_{\hat{X}^{s,\q,1}} ||\chi_{B_{\leq
\frac{i}{10}}} g^{K^c}||_{Y^{s}} \\ & \les \sum_{k} \sum_{i = k-2}^{k+2}
 2^{-\frac{i}{4}} ||f^K_{i}||_{\hat{X}^{s,\q,1}} ||
g||_{Y^{s}} \les ||f^K||_{\hat{X}^{s,\q,1}} ||
g^{K^c}||_{Y^{s}}.\end{align*} Finally, we use (\ref{b22}) to estimate
$I_2$:
\begin{align*}&||\frac{1}{\la \tau - |\xi|^2 \ra}I_{2}||_{X^{s,\q,1}} \les ||I_{2}||_{X^{s,-\q,1}} \\ &\les
\sum_{k} \sum_{d_3 \leq 2k-4} \sum_{i \geq k+11} \sum_{d_1 \leq
2i-4} 2^{ks} || (f_{i,d_{1}} * \frac{g^{K^c}}{w})_{k,d_3}||_{\hat{X}^{0,-\q}} \\
& \les \sum_{k} \sum_{d_3 \leq 2k-4} \sum_{i \geq k+11} \sum_{d_1
\leq 2i-4} 2^{ks}  2^{-4i}||f_{i,d_1}||_{\hat{X}^{0,\q}} ||
g^{K^c}||_{L^{2}}
\\ & \les
\sum_{k} \sum_{d_3 \leq 2k-4} \sum_{i \geq k+11} 2^{ks} 2^{-4i}
||f_{i}||_{\hat{X}^{0,\q,1}} || g^{K^c}||_{Y^s} \\ & \les \sum_{k} k
2^{ks} \sum_{i \geq k+11} 2^{-3i} ||f_{i}||_{\hat{X}^{s,\q,1}} ||
g^{K^c}||_{Y^s} \les  ||f||_{\hat{X}^{s,\q,1}} ||g^{K^c}||_{Y^s}. \end{align*}

Adding up the estimates we obtained for $I_{11}$, $I_{12}$ and
$I_{2}$ gives us the estimate (\ref{E3}).

\end{proof}

\section{The proof of the estimate \eqref{E4}}

We start with the following two elementary auxiliary lemmas.

\begin{lem} Let $f,g$ be non-negative smooth functions. Then,

\beq \label{b30} ||f_{i} * g_{j}||_{L^{2}_{\xi}L^{1}_{\tau}} \les
2^{\min{(i,j)}} ||f_{i}||_{L^{2}_{\xi}L^{1}_{\tau}} ||
g_{j}||_{L^{2}_{\xi}L^{1}_{\tau}}, \eeq

\beq \label{b31} ||(f_{i} * g_{j})_k||_{L^{2}_{\xi}L^{1}_{\tau}}
\les 2^{k} ||f_{i}||_{L^{2}_{\xi}L^{1}_{\tau}} ||
g_{j}||_{L^{2}_{\xi}L^{1}_{\tau}}. \eeq

Moreover if $|i-j| \leq 3$, $k \leq \max{(i,j)}-10$, and $d \leq
2\max{(i,j)}-10$ then

\beq \label{b32} ||(f^K_{i} *
\frac{g_{j}}{w})_{k,d}||_{L^{2}_{\xi}L^{1}_{\tau}} \les 2^{k}
2^{-20j} ||f_{i}||_{L^{2}_{\xi}L^{1}_{\tau}} ||
g_{j}||_{L^{2}_{\xi}L^{1}_{\tau}}. \eeq

\end{lem}

\begin{proof}

The first estimate is trivial

$$||f_{i} * g_{j}||_{L^{2}_{\xi}L^{1}_{\tau}} \les ||f_{i}||_{L^{1}_{\xi}L^{1}_{\tau}}
|| g_{j}||_{L^{2}_{\xi}L^{1}_{\tau}} \les  2^{i}
||f_{i}||_{L^{2}_{\xi}L^{1}_{\tau}} ||
g_{j}||_{L^{2}_{\xi}L^{1}_{\tau}}.$$

\noindent Analogously,  $$||f_{i} *
g_{j}||_{L^{2}_{\xi}L^{1}_{\tau}} \les  2^{j}
||f_{i}||_{L^{2}_{\xi}L^{1}_{\tau}} ||
g_{j}||_{L^{2}_{\xi}L^{1}_{\tau}},$$ and (\ref{b30}) immediately
follows.

For \eqref{b31} we estimate as follows:

\[
||(f_{i} * g_{j})_k||_{L^{2}_{\xi}L^{1}_{\tau}} \les 2^{k} ||(f_{i} * g_{j})_k||_{L^{\infty}_{\xi}L^{1}_{\tau}}
\les 2^{k} ||f_{i}||_{L^{2}_{\xi}L^{1}_{\tau}} ||
g_{j}||_{L^{2}_{\xi}L^{1}_{\tau}}.
\]

In order to derive (\ref{b32}) the key observation is that the
part of $g_j$ which interacts non-trivially is localized in a
region where $\tau \leq -2^{2j-10}$. Roughly speaking, we have a
high-high to low type of interaction, and $f^K_{i}$ is localized
at $\tau \geq 2^{2i-5}$. Using this observation, the argument to
prove \eqref{b32} is similar to the one for \eqref{b31}.

\end{proof}
\begin{lem} Let $f,g$ be non-negative smooth functions. Then,
\beq \label{b40} ||f_{i} * g_j||_{L^{2}} \les 2^{\min(i,j)}
||f_{i}||_{L^{2}_{\xi}L^{1}_{\tau}} || g_j||_{L^{2}}. \eeq

Moreover if $d \leq 2i-20$ then

\beq \label{b42} ||(1-\chi_{K}) (f^K_{i} *
\frac{g}{w})_d||_{L^{2}} \leq 2^{\frac{d}{2}} 2^{-20i}
||f_{i}||_{L^{2}_{\xi}L^{1}_{\tau}} || g||_{L^{2}}. \eeq

\end{lem}

\begin{proof}

The first estimate is trivial. Indeed,

$$||f_{i} * g_j||_{L^{2}} \les ||f_{i}||_{L^{1}_{\xi}L^{1}_{\tau}}
|| g_j||_{L^{2}} \les  2^{i} ||f_{i}||_{L^{2}_{\xi}L^{1}_{\tau}}
|| g_j||_{L^{2}},$$ and similarly we get the one with $2^j$
replacing $2^i$.

In order to derive (\ref{b42}), it is enough to observe that the
part of $g$ which gives nontrivial interactions is localized in a
region where $\tau \les -2^{2i-5}$. Moreover, the outcome of the
interaction is localized at $|\xi| \les 2^{\frac{d}{2}}$. These
two facts together with a similar argument to the one in the
previous Lemma give us the desired estimate.

To justify our first observation, we notice that $f^K_{i}$ is
localized at $\tau \geq 2^{2i-5}$, while $\chi_{B_d}(1-\chi_{K})$
is supported in a region where $|\tau| \leq 2^{6} |\tau-|\xi|^{2}|
\leq 2^{d+7} \leq 2^{2i-10}$. As a consequence the interacting
part of $g$ is supported in a region with $\tau \les -2^{2i-5}$.

To justify our second observation, we notice that the size of the
support of $\chi_{B_d}(1-\chi_{K})$ in the $\xi$ direction is
comparable to $2^{\frac{d}{2}}$, since on one hand we localize
$\langle \tau-|\xi|^2 \rangle$ around $2^{d}$, and on the other
hand we localize in a region where $|\tau-|\xi|^2| \geq
\frac{|\xi|^2}{2^{6}}$.

\end{proof}

\begin{proof}[Proof of $\eqref{E4}$] We wish to prove the
following estimate, \beq \label{E4new}
 \| \frac{(1-\chi_{K}) w}{\langle \tau - |\xi|^2 \rangle} (\frac{f^K}{w} * \frac{g^{K^c}}{w}) \|_{Y^{s}}
\lesssim \|f^K\|_{\hat{X}^{s,\q,1}} \|g^{K^c}\|_{Y^{s}}. \eeq

Since $w=1$ in $K$ we drop $w$ from the ratio $f^K/w$. We split
the proof in two steps.

\

\noindent \textbf{Step 1}. Estimate in $\la \xi \ra^{s}
L^{2}_{\xi}L^{1}_{\tau}.$

In this proof we use the $\hat{X}^{s,\q,1}$ structure for $f^K$
via the inclusion $X^{s,\q,1} \subset \la \xi \ra^{s} L^{2}_{\xi}
L^{1}_{\tau}$. We decompose
\begin{align} \label{dec10} &(1-\chi_{K})(f^K * \frac{g^{K^{c}}}{w})= \sum_{k}
\sum_{i} \sum_{j} \chi_{K^{c}} (f^K_{i} *
\frac{g^{K^{c}}_{j}}{w})_k\\& = \sum_{i \leq j} \sum_{k =
j-10}^{j+2} \chi_{K^{c} } (f^K_{i} * \frac{g^{K^{c}}_{j}}{w})_k +
\mbox{symmetric term} \nonumber \\ & +\sum_{i} \sum_{j} \sum_{k
\leq \max{(i,j)-10}}  \chi_{K^{c} } (f^K_{i} *
\frac{g^{K^{c}}_{j}}{w})_k\nonumber
\\ &= J_1 + \mbox{symmetric term} + J_2.\nonumber\end{align}
When estimating $J_1$ we drop all the weights and use (\ref{b30})
to obtain
\begin{align*} & ||\frac{w J_1}{\la \tau - |\xi|^2 \ra} ||^2_{\la \xi \ra^{s} L^{2}_{\xi}L^{1}_{\tau}}
\les  \sum_{j} \sum_{k = j-10}^{j+2} 2^{2ks} 2^{-4k} ||
\sum_{i=0}^j \chi_{K^{c}} (f^K_{i} * g_{j}^{K^c})_k||^2_{
L^{2}_{\xi}L^{1}_{\tau}} \\ & \les \sum_{j} \sum_{k = j-10}^{j+2}
2^{-4k} 2^{2(k-j)s} \left ( \sum_{i=0}^j 2^{i(1-s)}
||f^K_{i}||_{\la \xi \ra^{s} L^{2}_{\xi}L^{1}_{\tau}} \right)^2
||g_{j}^{K^c}||^2_{\la \xi \ra^{s} L^{2}_{\xi}L^{1}_{\tau}} \\ &
\les \sum_{j} \sum_{k = j-10}^{j+2} \left( \sum_{i=0}^j
2^{i(1-s)-2j} ||f^K_{i}||_{\la \xi \ra^{s}
L^{2}_{\xi}L^{1}_{\tau}} \right)^2 ||g_{j}^{K^c}||^2_{\la \xi
\ra^{s} L^{2}_{\xi}L^{1}_{\tau}} \\ & \les\sum_{j} \sum_{k =
j-10}^{j+2} ||f^{K}||^2_{\la \xi \ra^{s} L^{2}_{\xi}L^{1}_{\tau}}
||g_{j}^{K^c}||^2_{\la \xi \ra^{s} L^{2}_{\xi}L^{1}_{\tau}} \les
||f^{K}||_{\la \xi \ra^{s} L^{2}_{\xi}L^{1}_{\tau}}
||g^{K^c}||_{\la \xi \ra^{s} L^{2}_{\xi}L^{1}_{\tau}}.\end{align*}
For $J_2$, we notice that in order to consider only the nontrivial
terms, we must impose the condition $|i-j| \leq 3$. We further
split
\begin{align*}J_2 &=  \sum_{|i-j| \leq 3} \sum_{k \leq \max{(i,j)-10}} \sum_{d \geq 2\max{(i,j)-10}}
\chi_{K^{c}} (f^K_{i} * \frac{g_{j}^{K^c}}{w})_{k,d}
\\&+\sum_{|i-j| \leq 3} \sum_{k \leq \max{(i,j)-10}} \sum_{d \leq
2\max{(i,j)-10}} \chi_{K^{c}} (f^K_{i} *
\frac{g_{j}^{K^c}}{w})_{k,d}\\&=J_{21} + J_{22}.\end{align*}The
term  $J_{21}$ is estimated using (\ref{b31}):\begin{align*} &
||\frac{w J_{21}}{\la \tau - |\xi|^2 \ra}||_{\la \xi \ra^{s}
L^{2}_{\xi}L^{1}_{\tau}} \les \sum_{|i-j| \leq 3} \sum_{k \leq
\max{(i,j)-10}} 2^{ks} 2^{-2 j} || \chi_{K^{c}} (f^K_{i}
* g_{j}^{K^c})_k||_{\la \xi \ra^{s} L^{2}_{\xi}L^{1}_{\tau}} \\ &\les
\sum_{|i-j| \leq 3} \sum_{k \leq \max{(i,j)-10}} 2^{(s+1)k}
2^{-2(s+1)j} ||f^K_{i}||_{\la \xi \ra^{s} L^{2}_{\xi}L^{1}_{\tau}}
||g_{j}^{K^c}||_{\la \xi \ra^{s} L^{2}_{\xi}L^{1}_{\tau}} \\ &\les
\sum_{|i-j| \leq 3} ||f^K_{i}||_{\la \xi \ra^{s}
L^{2}_{\xi}L^{1}_{\tau}} ||g_{j}^{K^c}||_{\la \xi \ra^{s}
L^{2}_{\xi}L^{1}_{\tau}} \les ||f^K||_{\la \xi \ra^{s}
L^{2}_{\xi}L^{1}_{\tau}} ||g^{K^c}||_{\la \xi \ra^{s}
L^{2}_{\xi}L^{1}_{\tau}}.\end{align*} The term $J_{22}$ is bounded
using (\ref{b32}) and keeping the weights,
\begin{align*} & ||\frac{w}{\la \tau - |\xi|^2 \ra} \sum_{|i-j| \leq 3} \sum_{k \leq \max{(i,j)-10}}
\sum_{d \leq 2\max{(i,j)-10}} \chi_{K^{c}} (f^K_{i} *
\frac{g_{j}^{K^c}}{w})_{k,d}||_{\la \xi \ra^{s}
L^{2}_{\xi}L^{1}_{\tau}}
\\&\les\sum_{|i-j| \leq 3} \sum_{k \leq \max{(i,j)-10}}
\sum_{d \leq 2\max{(i,j)-10}}  2^{ks} 2^{9d} || (f^K_{i} *
\frac{g_{j}^{K^c}}{w})_{k,d}||_{L^{2}_{\xi}L^{1}_{\tau}} \\&\les
\sum_{|i-j| \leq 3} \sum_{k \leq \max{(i,j)-10}} \sum_{d \leq
2\max{(i,j)-10}} 2^{(s+1)k} 2^{9d} 2^{-2js}
2^{-20j}||f^K_{i}||_{\la \xi \ra^{s} L^{2}_{\xi}L^{1}_{\tau}}
||g_{j}^{K^c}||_{\la \xi \ra^{s} L^{2}_{\xi}L^{1}_{\tau}}
\\&\les\sum_{|i-j| \leq 3} \sum_{k \leq \max{(i,j)-10}}
2^{(s+1)k} 2^{-2(s+1)j} ||f^K_{i}||_{\la \xi \ra^{s}
L^{2}_{\xi}L^{1}_{\tau}} ||g_{j}^{K^c}||_{\la \xi \ra^{s}
L^{2}_{\xi}L^{1}_{\tau}} \\&\les \sum_{|i-j| \leq 3}
||f^K_{i}||_{\la \xi \ra^{s} L^{2}_{\xi}L^{1}_{\tau}}
||g_{j}^{K^c}||_{\la \xi \ra^{s} L^{2}_{\xi}L^{1}_{\tau}} \les
||f^K||_{\la \xi \ra^{s} L^{2}_{\xi}L^{1}_{\tau}} ||g^{K^c}||_{\la
\xi \ra^{s} L^{2}_{\xi}L^{1}_{\tau}}.\end{align*}

Combining the estimates for $J_1$, $J_{21}$ and $J_{22}$ we obtain
the claim in (\ref{E3}).

\

\noindent \textbf{Step2}. Estimate in $\la (\xi,\tau) \ra^{s+1}
L^{2}.$

We decompose
\begin{align*} & (1-\chi_{K})(f^K * \frac{g^{K^c}}{w})= \sum_{d_3} \sum_{i} \chi_{K^{c} } (f^K_{i} * \frac{g^{K^c}}{w})_{d_3} \\ &=
\sum_{i} \sum_{d_3 \geq 2i-10}  \chi_{K^{c} } (f^K_{i}
* \frac{g^{K^c}}{w})_{d_3}+ \sum_{i} \sum_{d_3 \leq 2i-10}
\chi_{K^{c}} (f^K_{i} * \frac{g^{K^c}}{w})_{d_3}\\&=L_1 +
L_2.\end{align*}

We estimate $L_1$ by using (\ref{b40}) and dropping all the
weights:
\begin{align*} & ||\frac{w L_1}{\la \tau - |\xi|^2 \ra}||_{\la (\xi, \tau) \ra^{1+s}  L^{2}} \les
\sum_{i} \sum_j \sum_{d_3 \geq 2i-10} 2^{s d_3} || (f^K_{i} *
g_j^{K^c})_{d_3}||_{L^{2}} \\&\les \sum_{i} \sum_j 2^{si}
2^{\min{(i,j)}} 2^{-(1+s)j} \sum_{d_3 \geq 2i-10} 2^{s (d_3-2i)}
||f_{i}^{K}||_{\la \xi \ra^{s} L^{2}_{\xi}L^{1}_{\tau}}
||g_j^{K^c}||_{\la (\xi, \tau) \ra^{1+s}L^{2}} \\&\les \sum_{i
\leq j} 2^{(s+1)(i-j)} ||f_{i}^{K}||_{\la \xi \ra^{s}
L^{2}_{\xi}L^{1}_{\tau}} ||g_j^{K^c}||_{\la (\xi, \tau)
\ra^{1+s}L^{2}}
\\&+ \sum_{i \geq j} 2^{s(i-j)} ||f_{i}^{K}||_{\la \xi \ra^{s}
L^{2}_{\xi}L^{1}_{\tau}} ||g_j^{K^c}||_{\la (\xi, \tau)
\ra^{1+s}L^{2}}\\& \les||f^{K}||_{\la \xi \ra^{s}
L^{2}_{\xi}L^{1}_{\tau}} ||g^{K^c}||_{\la (\xi, \tau)
\ra^{1+s}L^{2}}.\end{align*}

In order to estimate $L_2$ we use (\ref{b42}) and the fact that
$g$ is supported where $\tau \leq -2^{2i-5}$. Hence,
\begin{align*} & ||\frac{w(1-\chi_{K})}{\la \tau - |\xi|^2 \ra} \sum_{i}
\sum_{d_3 \leq 2i-10} (f^K_{i} * \frac{g}{w})_{d_3}||^{2}_{\la
(\xi, \tau) \ra^{1+s}  L^{2}}\\& \les||\frac{w(1-\chi_{K})}{\la
\tau - |\xi|^2 \ra} \sum_{d_3} \sum_{i \geq \frac{d_3+10}{2} }
(f^K_{i} * \frac{g}{w})_{d_3}||^{2}_{\la (\xi, \tau) \ra^{1+s}
L^{2}}
\\& \les\sum_{d_3} \sum_{i \geq \frac{d_3+10}{2}} 2^{(s+19) d_3}
|| \chi_{K^{c}} (f^K_{i} * \frac{g}{w})_{d_3}||^2_{L^{2}}
\\& \les \sum_{d_3} \sum_{i \geq \frac{d_3+10}{2}} 2^{-2is} 2^{(s+19) d_3}
2^{d_3} 2^{-40i} 2^{-2(s+1)i}||f_{i}||^2_{\la \xi \ra^{s}
L^{2}_{\xi}L^{1}_{\tau}} ||g||^2_{\la (\xi, \tau) \ra^{1+s}L^{2}}
\\& \les \sum_{d_3} \sum_{i \geq \frac{d_3+10}{2}} 2^{(s+20)
(d_3-2i)} 2^{-2i(s+1)} ||f_{i}||^{2}_{\la \xi \ra^{s}
L^{2}_{\xi}L^{1}_{\tau}} ||g||^2_{\la (\xi, \tau)
\ra^{1+s}L^{2}}\\& \les ||f||^2_{\la \xi \ra^{s}
L^{2}_{\xi}L^{1}_{\tau}} ||g||^2_{\la (\xi, \tau)
\ra^{1+s}L^{2}},\end{align*} which concludes our proof.
\end{proof}

\end{document}